\documentclass[11pt]{amsart}
\usepackage{amsmath}
\usepackage[all]{xy}
\usepackage{amssymb}
\usepackage{mathrsfs}
\usepackage[noadjust]{cite}
\usepackage{hyperref}
\usepackage{amsfonts}
\usepackage{mathdots}
\usepackage{todonotes}
\usepackage{bm}
\usepackage{tikz}
\newcommand*\circled[1]{\tikz[baseline=(char.base)]{
            \node[shape=circle,draw,inner sep=2pt] (char) {#1};}}

\theoremstyle{plain}
\newtheorem{thm}{Theorem}[section]
\newtheorem{prop}[thm]{Proposition}
\newtheorem{lem}[thm]{Lemma}
\newtheorem{cor}[thm]{Corollary}

\theoremstyle{remark}
\newtheorem{rem}[thm]{Remark}

\theoremstyle{definition}
\newtheorem{defn}[thm]{Definition}

\newtheorem{Ques}[thm]{Question}
\newcommand{\cA}{\mathcal{A}}
\theoremstyle{conjecture}

\newcommand{\Ku}{\mathcal{K}u}

\def\cO{\mathcal{O}}

\def\HH{\mathrm{HH}}
\def\Hom{\mathrm{Hom}}
\def\Id{\mathrm{Id}}
\def\Jac{\mathrm{Jac}}

\newcommand{\comments}[1]{}

\newcommand{\into}{\hookrightarrow}
\newcommand{\bR}{\bm{\mathrm{R}}}
\newcommand{\bL}{\bm{\mathrm{L}}}
\renewcommand{\AA}{\ensuremath{\mathbb{A}}}
\newcommand{\CC}{\ensuremath{\mathbb{C}}}

\newcommand{\PP}{\ensuremath{\mathbb{P}}}

\newcommand{\ZZ}{\ensuremath{\mathbb{Z}}}
\newcommand{\cB}{\ensuremath{\mathcal{B}}}
\newcommand{\cC}{\ensuremath{\mathcal{C}}}

\newcommand{\cE}{\ensuremath{\mathcal{E}}}

\DeclareMathOperator{\id}{\mathrm{id}}
\DeclareMathOperator{\rk}{\mathrm{rk}}

\DeclareMathOperator{\GL}{GL}

\DeclareMathOperator{\Fact}{Fact}
\DeclareMathOperator{\Acyc}{Acyc}
\DeclareMathOperator{\Inj}{Inj}

\newcommand{\abs}{\text{abs}}

\DeclareMathOperator\oh{\mathcal{O}}
\date{}

\title[Categorical Torelli theorems for weighted hypersurfaces]{IVHS via Kuznetsov components and categorical Torelli theorems for weighted hypersurfaces}

\date{}

\begin{document}

\author{Xun Lin, Jørgen Vold Rennemo and Shizhuo Zhang}

\begin{abstract}
We study the categorical Torelli theorem for smooth (weighted) hypersurfaces in (weighted) projective spaces via the Hochschild--Serre algebra of its Kuznetsov component.
In the first part of the paper, we show that a natural graded subalgebra of the Hochschild--Serre algebra of the Kuznetsov component of a degree $d$ weighted hypersurface in $\mathbb{P}(a_0,\ldots,a_n)$ reconstructs the graded subalgebra of the Jacobian ring generated by the degree $t:=\mathrm{gcd}(d,\Sigma_{i=0}^na_i)$ piece under mild assumptions. 
Using results of Donagi and Cox--Green, this gives a categorical Torelli theorem for most smooth hypersurfaces $Y$ of degree $d \le n$ in $\mathbb{P}^n$ such that $d$ does not divide $n+1$ (the exception being the cases of the form $(d,n) = (4, 4k + 2)$, for which a result of Voisin lets us deduce a generic categorical Torelli theorem when $k \ge 150$).
Next, we show that the Jacobian ring of the Veronese double cone can be reconstructed from its graded subalgebra of even degree, thus proving a categorical Torelli theorem for the Veronese double cone. 

In the second part, we rebuild the infinitesimal Variation of Hodge structures of a series of (weighted) hypersurfaces from their Kuznetsov components via the Hochschild--Serre algebra. As a result, we prove categorical Torelli theorems for two classes of (weighted) hypersurfaces: $(1):$ Generalized Veronese double cone; $(2):$ Certain $k$-sheeted covering of $\mathbb{P}^n$, when they are generic. Then, we prove a refined categorical Torelli theorem for a Fano variety whose Kuznetsov component is a Calabi--Yau category of dimension $2m+1$. Finally, we prove the actual categorical Torelli theorem for generalized Veronese double cone and $k$-sheeted covering of $\mathbb{P}^n$.
\end{abstract}

\address{Max Planck Institute for Mathematics, Vivatsgasse 7, 53111 Bonn, Germany}
\email{xlin@mpim-bonn.mpg.de}

\address{Department of Mathematics, University
of Oslo, PO Box 1053 Blindern, Oslo 0316,
Norway}
\email{jorgeren@uio.no}

\address{Simons Laufer Mathematical Sciences Institute,17 Gauss Way, Berkeley, CA 94720}
\address{Morningside Center of Mathematics, Chinese Academy of Science}
\email{shizhuozhang@msri.org, shizhuozhang@mpim-bonn.mpg.de}


\maketitle

\section{introduction}

Homological algebra is both a fundamental tool for studying sheaves on algebraic varieties and a source of interesting varieties' invariants. 
The bounded derived category of coherent sheaves on an algebraic variety can be viewed as a universal homological invariant of an algebraic variety. Thus, a natural question is whether the bounded derived category of coherent sheaves on an algebraic variety determines its isomorphism class. The seminal work \cite{bondal2001reconstruction} established a reconstruction theorem for smooth Fano and general type variety via its derived category, and this is the program's starting point. 

Let $X$ be a smooth projective variety.
The \emph{Torelli problem} asks if $X$ is determined up to isomorphism (within its deformation equivalence class) by a Hodge structure attached to $X$.
As $D^{b}(X)$ is a very rich invariant, it is not surprising that in certain cases the Hodge theory of $X$ can be reconstructed from the derived category $D^{b}(X)$.\footnote{But note that e.g.~if $X$ is a K3 surface, the integral Hodge structure $H^2(X)$ is not determined by $D^{b}(X)$ \cite{orlov1997equivalences}.}
We denote by ``$\xymatrix{A\ar[r]^{?}&B\\}$"  the question ``Does $A$ reconstruct $B$?". Then, we can ask the following question for a given deformation class of varieties.


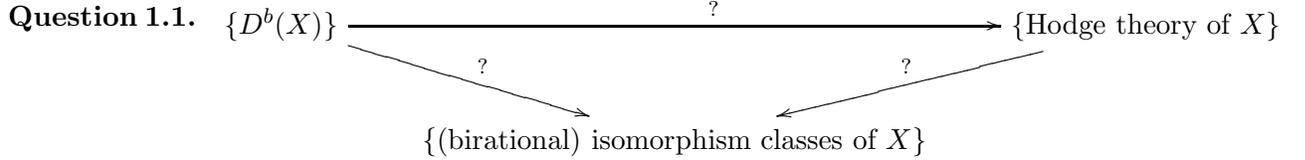
\begin{Ques}\label{commutativediagram1}\ 
\xymatrix{\{\text{$D^{b}(X)$}\}\ar[rr]^{?}\ar[dr]^{?}&&\{\text{Hodge theory of $X$}\}\ar[dl]_{?}\\
&\{\text{(birational) isomorphism classes of $X$}\}&}
\end{Ques}
In the case where $X$ is a Fano variety, a choice of an exceptional collection of objects on $X$ gives rise to an admissible subcategory $\Ku(X)\subset D^{b}(X)$, defined as the left orthogonal to the exceptional collection.
This category is known as the \emph{Kuznetsov component}, and often encodes essential information about $X$. 
Then we are interested in the following variation of Question~\ref{commutativediagram1}.

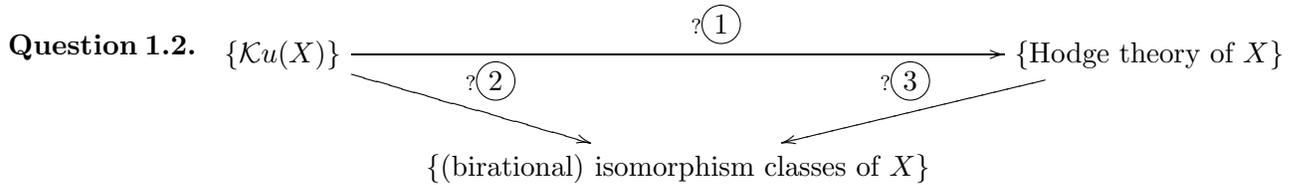
\begin{Ques}\label{Ques}\ 
\xymatrix{\{\text{$\Ku(X)$}\}\ar[rr]^{? 
 \circled{1}}\ar[dr]^{? \circled{2}}&&\{\text{Hodge theory of $X$}\}\ar[dl]_{? \circled{3} }\\
&\{\text{(birational) isomorphism classes of $X$}\}&}
\end{Ques}

Here question \circled{3} is known as the \emph{(birational) Torelli problem}. 
Tremendous work is carried out to establish the Torelli theorem for a large class of Fano varieties (see \cite{CG72,debarre1989theoreme, debarre1990theoreme,Donagi1986,voisin2022schiffer}). 
Question \circled{2} is known as the \emph{(birational) categorical Torelli problem}. 
There are several deformation classes of smooth Fano varieties for which the answer to question \circled{2} is known to be positive, such as the cubic threefold (\cite{bernardara2012categorical}) and the intersection of two quadrics in $\mathbb{P}^5$ (\cite{kuznetsov2009derived}). 
The paper \cite{pertusi2023categorical} gives an overview of recent results in this direction. 
As the diagram illustrates, question \circled{2} may be given a positive answer by giving suitable positive answers to questions \circled{1} and \circled{3}. 
There are several ways to give positive answers to question \circled{1}. 
In \cite{bernardara2013derived} and \cite{perry2022integral}, the authors reconstruct the intermediate Jacobian of an odd-dimensional smooth Fano variety from its Kuznetsov component. 
In \cite{Dell2023categorical}, the authors reconstruct the primitive cohomology of a smooth projective surface from the equivariant Kuznetsov component. 
In \cite{lin2024serre}, the first and third author of this paper reconstruct the Jacobian ring of a series of Fano hypersurfaces in $\mathbb{P}^n$ and deduce a \emph{Categorical Torelli theorem} for them. 
In our paper, we continue to study question \circled{1} and introduce a new technique to reconstruct the infinitesimal variation of Hodge structures(IVHS) of a series (weighted) hypersurfaces from their Kuznetsov components and deduce \emph{categorical Torelli theorem} for them. 

 For a smooth Fano hypersurface $X\subset\mathbb{P}^n$ of degree $d$, the generic Torelli theorem has been established except for finitely many values of $d$ and $n$ (\cite{Don, CM_1990__73_2_121_0, voisin2022schiffer}). 
 In \cite{Don}, a key idea is to recover a sub-graded ring of the Jacobian ring $\mathrm{Jac}(f)$ from the infinitesimal variation of Hodge structure of $X$ using certain ``symmetrizer lemma'' \cite[Lemma 6.3]{Don}. 
 Next, Donagi shows that knowledge of this sub-graded ring is enough to determine the surjection $\CC[x_0,\dots, x_n] \to \mathrm{Jac}(f)$ up to $\mathrm{GL}(n+1)$-equivalence.
 By the Mather--Yau theorem, this suffices to determine $X$ up to projective equivalence.
 Since the infinitesimal variation of Hodge structure of $X$ encodes the differential of the map from the moduli space of hypersurfaces to the period domain, it follows by the ``principle of prolongation'' that the generic Torelli theorem holds.
 
 In the case of hypersurfaces in \textit{weighted} projective spaces, a weighted symmetrizer lemma is proved in \cite[Theorem 1.6]{Donagi1986}. Still, the necessary hypotheses in this version means that there are many interesting examples for which the lemma does not apply.
 
 In this paper, we show and exploit the fact that the Kuznetsov component $\Ku(X)$ reconstructs the IVHS of $X$ for many interesting examples. 
 Starting from $\Ku(X)$, there is no need for a symmetrizer lemma, and the relevant sub-graded ring of the Jacobian ring can instead be extracted from the Hochschild--Serre algebra of $\Ku(X)$.
 We consider the following examples.

\begin{enumerate}
\item Smooth projective hypersurfaces in $\PP^n$.
    \item\label{enum:veroneseDoubleCone} Veronese double cones (i.e.,~smooth del Pezzo threefolds of Picard number one and degree one).
    \item Generalized Veronese double cones (see Section~\ref{section_GVDC}).
    \item Certain $k$-sheeted coverings of $\mathbb{P}^n$ (see Section \ref{section_kSC}).
\end{enumerate}
The threefold in (\ref*{enum:veroneseDoubleCone}) can be described as a degree 6 hypersurface in $\mathbb{P}(1,1,1,2,3)$.

\subsection{Smooth hypersurfaces in projective space}
Let $Y$ be a hypersurface of $\mathbb{P}^n$ of degree $d\leq n$. 
In this case the Kuznetsov component is given by 
 $$\Ku(Y):=\Big\langle \mathcal{O}_{Y},\mathcal{O}_{Y}(1),\cdots, \mathcal{O}_{Y}\left(n-d\right)\Big\rangle^{\perp}.$$
The first main result of our article is \emph{Categorical Torelli theorem} for them. 

\begin{thm}\label{theorem_first_main_result_categorical_Torelli_hypersurfaces}(=Theorem~\ref{thm_rennemo2023}+Corollary~\ref{cor_generic_categorical_Torelli_quartic_hypesurface})
\leavevmode\begin{enumerate}
    \item Let $X$ and $X'$ be smooth hypersurfaces of degree $d \le n$ in $\mathbb{P}^{n}$.
 Assume that $d$ does not divide $n+1$, and that $(d, n+1)$ is not of the form $(4,4k+2)$ for $k \in \ZZ$. 
 If $\Ku(X) \simeq\Ku(X')$ and the equivalence is given by a Fourier--Mukai functor, then $X\cong X'$.
 \item Let $X$ and $X'$ be generic smooth quartic hypersurfaces in $\mathbb{P}^{4k+1}$, with $k \ge 150$.
If there is a Fourier--Mukai type equivalence $\Phi:\Ku(X)\simeq\Ku(X')$, then $X\cong X'$. 
\end{enumerate}
\end{thm}

\subsection{Veronese double cone and Kuznetsov component}
A Veronese double cone $Y$ is a hypersurface cut out by a sextic equation in the weighted projective space $\mathbb{P}:=\mathbb{P}(1,1,1,2,3)$. Let $x_1,\ldots,x_5$ be coordinates of $\mathbb{P}$, where $x_4,x_5$ are those of weight $2$ and $3$ respectively. By completing a square, we can write the equation for $Y$ as $x_5^2=f(x_1,x_2,\ldots,x_4)$, where $f$ is a degree 6 polynomial. From \cite[Proposition 3.1]{iskovskih1977fano}, we know the linear system $|-K_Y|=|2H|$ is basepoint free and induces a morphism $\phi_{2H}\colon Y\rightarrow\mathbb{P}(H^0(\oh_Y(2)))=\mathbb{P}^6$, whose image $K\cong\mathbb{P}(1,1,1,2)$ is the cone over the Veronese surface. Indeed, the morphism $\phi_{2H}$ is smooth of degree two outside the vertex of the cone and the divisor $W:=\{f=0\}$, which justifies the name of the Veronese double cone. We write the equation of the branch divisor $W$ as $\omega=x_4^3+g(x_1,x_2,x_3)x_4+h(x_1,x_2,x_3)$, where degree of $g$ and $h$ are $4$ and $6$ respectively. The variety $Y$ has a semi-orthogonal decomposition given by 
$$D^b(Y)=\langle\Ku(Y),\oh_Y,\oh_Y(1)\rangle.$$ 
The second main result of our article is \emph{Categorical Torelli theorem} for Veronese double cone. 

\begin{thm}\label{thm_main_result_categorical_Torelli_theorem}(=Theorem~\ref{thm_categorical_Torelli_theorem_for_Veronese_double_cone})
Let $Y$ and $Y'$ be Veronese double cones such that $\Ku(Y)\simeq\Ku(Y')$ and the equivalence is given by a Fourier--Mukai functor, then $Y\cong Y'$. 
\end{thm}

A popular strategy for proving categorical Torelli theorems is to introduce a stability condition on $\Ku(Y)$, which in turn enables the construction of moduli spaces of stable objects, from which one can recover some specific classical moduli space.
Then, one can recognize the isomorphism class of the varieties via the classical moduli space. Indeed, by \cite{bayer2017stability}, the Kuznetsov component of 
any del Pezzo threefold of Picard rank 1 admits a Serre-invariant stability condition, and several interesting moduli spaces are constructed from $\Ku(Y)$ (see \cite{petkovic2023note}). But for higher dimensional smooth hypersurface $Y$ in projective space, it is not known whether stability conditions exist on $\Ku(Y)$. Even if the Serre-invariant stability condition exists on Kuznetsov component $\Ku(Y)$ of Veronese double cone $Y$, the orbit of Serre-invariant stability conditions is not unique, so the moduli space constructed from $\Ku(Y)$ is not canonical, thus it is difficult to relate the moduli space constructed from $\Ku(Y)$ to a classical moduli space over $Y$ and compare it to that constructed from $\Ku(Y')$ whenever there is an equivalence $\Ku(Y)\simeq\Ku(Y')$.  Another approach to solving the categorical Torelli problem is to reduce it to a classical Torelli problem in terms of intermediate Jacobians via the abstract construction of intermediate Jacobian of admissible subcategories of smooth projective varieties (see \cite{perry2022integral}). However, the lack of classical Torelli theorem for the Veronese double cone in the literature makes it impossible to work this way.


To overcome these difficulties, we introduce the new techniques of Hochschild--Serre algebra $\mathrm{HS}(\Ku_{dg}(Y))$ constructed from the dg enhancement $\Ku_{dg}(Y)$ of its Kuznetsov component when $Y$ is a smooth hypersurface in projective space or a weighted hypersurface (the Veronese double cone is regarded as a weighted hypersurface of degree $6$ in $\mathbb{P}(1,1,1,2,3)$). 
By \cite{BFK}, the Kuznetsov category is equivalent to a graded matrix factorization category, and moreover its Hochschild--Serre algebra is computed explicitly.
This makes it possible to relate an explicit sub-algebra of the Hochschild--Serre algebra to the Jacobian ring of the polynomial defining $Y$. 
The upshot is that, for the (weighted) hypersurfaces above, we can always recover some graded subalgebra of their Jacobian ring using the Hochschild--Serre algebra of their Kuznetsov components. In the first case, 
we apply  techniques in \cite{Don,CM_1990__73_2_121_0,voisin2022schiffer} to recover the surjection $\mathbb{C}[x_0,\ldots,x_n]\rightarrow\mathrm{Jac}(f)$ up to $\mathrm{GL}(n+1)$-equivalence. Then, we use Mather-Yau's theorem to determine the isomorphism class, thus Theorem~\ref{theorem_first_main_result_categorical_Torelli_hypersurfaces} is established. 
In the second case, we recover the even-degree subalgebra of the Jacobian algebra, and apply additional ad hoc techniques to extend the reconstruction to the whole Jacobian ring. Thus Theorem~\ref{thm_main_result_categorical_Torelli_theorem} is established.

\begin{rem}\label{rem_Franco}
In \cite{Dell2023categorical}, the authors also proved a version of the categorical Torelli theorem for Veronese double cone with the additional assumption that the equivalence of Kuznetsov components commutes with rotation functors and Veronese double cones are very general. 
\end{rem}

\subsection{Reconstruction of IVHS from Kuznetsov components}\label{section_reconstruction_IVHS}
Let $X$ be a smooth projective variety over $\mathbb{C}$. The so-called Infinitesimal Variation of Hodge structures (IVHS) \cite{carlson2017period} consists of integral polarised Hodge structure as transcendental data and 
$$\{H^1(X,T_X),\bigoplus_{p+q=n}\mathrm{Hom}(H^p(X,\Omega_X^q),H^{p+1}(X,\Omega_X^{q-1})), d\mathcal{P}\}$$ as algebraic data. By \cite{carlson2017period} the differential $d\mathcal{P}$ of the period map is given by 
$$ d\mathcal{P} : H^1(X,T_X) \longrightarrow \bigoplus_{p+q=n}\mathrm{Hom}(H^p(X,\Omega_X^q),H^{p+1}(X,\Omega_X^{q-1})).$$

When we say IVHS, we always refer to its algebraic data.

\subsubsection{$\mathrm{IVHS}$ for (weighted) hypersurfaces}
In the case of a (weighted) hypersurface $Y\subset\mathbb{P}(a_0,\ldots,a_n)$ with $\mathrm{gcd}(a_0,\ldots,a_n)=1$ and dimension greater than two, the IVHS for $Y$ has a particularly nice form: it is shown in Lemma~\ref{weightedIVHS} that the IVHS of $Y$ is the graded component of its Jacobian ring $\mathrm{Jac}(\omega)$. The importance of the IVHS is that one can use it to determine the isomorphism class of those (weighted) hypersurfaces at least generically, as is shown in \cite{Donagi1986},\cite{CM_1990__73_2_121_0}, \cite{Don} and \cite{voisin2022schiffer}. The third main result of our article is a reconstruction of IVHS of a class of (weighted) hypersurfaces via their Kuznetsov components. 



\begin{thm}\label{thm_second_main_result}(=Proposition~\ref{prop_reconstruction_IVHS_quartic_hypersurface}+Theorem~\ref{veroneseIVHS}+Theorem~\ref{thm_generic_Torelli_2}).
The infinitesimal Variation of Hodge structures (IVHS) for the following (weighted) hypersurfaces is reconstructed from their Kuznetsov components. More precisely, if the equivalence $\Phi\colon \Ku(Y)\simeq\Ku(Y')$ is of Fourier--Mukai type, then $\mathrm{IVHS}(Y)\cong\mathrm{IVHS}(Y')$\footnote{The isomorphism of IVHS is the isomorphism of graded pieces of the Jacobian rings, which is compatible with the graded multiplication.}
\begin{enumerate}
\item $Y\subset\mathbb{P}^n$ of degree $3\leq d\leq n$ and $(d,n+1)\neq (3,6)$. 
    \item $Y\subset\mathbb{P}(\underbrace{1,\ldots,1}_{r+1},2n,kn)$ with degree $2kn$ such that $k\geq 3$ is odd and $r\geq 2$ and $m|2n$, where $m=\mathrm{gcd}((k+2)n+r+1,2kn)$.
    \item $X\subset\mathbb{P}(1,1,\cdots,1,s)$ of degree $ks$ and $r,s,k\geq 2$ with $(r,s,k)$ is not $(2,2,2)$ and $$ 2\cdot n_{0}<\operatorname
{min}\left(2s, ((r+2)k-2)s-2(r+1),(ks-2)(r+1)\right),$$ 
and $n_{0}|s$, where $n_{0}=gcd(r+1+s,ks)$. 

\end{enumerate}
\end{thm}

As a corollary, we prove the categorical Torelli theorem for generic (weighted) hypersurfaces. 
\begin{cor}\label{thm_the_third_main_result}(Theorem~\ref{generictorelli}+Theorem~\ref{thm_generic_Torelli_2}).
Categorical Torelli theorems hold for following (weighted) hypersurfaces:\begin{enumerate}
    \item Degree $2kn$ generic hypersurface $Y\subset\mathbb{P}(\underbrace{1,\ldots,1}_{r+1},2n,kn)$ satisfying assumptions in Theorem~\ref{thm_second_main_result} $(1)$ and $m<2n$. 
    \item Degree $ks$ generic hypersurface in $\mathbb{P}(1,1,\cdots,1,s)$ satisfying assumptions in Theorem~\ref{thm_second_main_result} $(3)$. 
\end{enumerate}
More precisely, if $Y$ and $Y'$ are one of those (weighted) hypersurfaces above and $\Phi\colon\Ku(Y)\simeq\Ku(Y')$ is an equivalence of Fourier--Mukai type, then $Y\cong Y'$. 
\end{cor}



We call a Fano variety whose Kuznetsov component is a Calabi--Yau category of dimension $n$ a Fano variety of $\mathrm{CYn}$-type. For example, cubic fourfolds and Gushel-Mukai fourfolds are Fano varieties of $\mathrm{CY2}$-type, and a cubic sevenfold is of $\mathrm{CY3}$-type. In \cite{huybrechts2016hochschild} and \cite{bayer2017stability}, two different proofs of a refined version of categorical Torelli theorem for cubic fourfolds are given, i.e. the Kuznetsov component of a cubic fourfold, together with its rotation functor determines its isomorphism class. Next we consider a double covering of $\mathbb{P}^{4n+1}$, branched over a generic quartic hypersurface. Its Kuznetsov component is a $\mathrm{CY-(2n+1)}$ category. We prove a refined categorical Torelli theorem for them:

\begin{thm}\label{theorem_refined_categorical_Torelli}(=Corollary~\ref{cor_categorical_Torelli_double_cover})
Let $Y,Y'\rightarrow\mathbb{P}^{4n+1}$ be the double covers of $\mathbb{P}^{4n+1}$ branched over generic quartic hypersurfaces such that $n\geq 150$. Assume that there is a Fourier--Mukai type equivalence $\Phi:\Ku(Y)\simeq\Ku(Y')$ such that $\Phi$ commutes with rotation functors $\bR$ and $\bR'$ respectively. Then $Y\cong Y'$. 
\end{thm}

\subsection{Categorical Torelli Theorem for two series of weighted hypersurfaces}
we prove the categorical Torelli theorem for weighted hypersurfaces in Theorem ~\ref{thm_second_main_result} if they are generic. In this section, we state an improvement of the results by directly reconstructing the Jacobian ring of these weighted hypersurfaces from their Kuznetsov components.

\begin{thm}\label{thm_actual_CTT}(=Theorem~\ref{thm_categorical_Torelli_thm_weighted_hypersurfaces})
Categorical Torelli theorems hold for following (weighted) hypersurfaces:\begin{enumerate}
    \item Degree $2kn$ hypersurface $Y\subset\mathbb{P}(\underbrace{1,\ldots,1}_{r+1},2n,kn)$ satisfying assumptions in Theorem~\ref{thm_second_main_result} $(1)$ and $m<2n$. 
    \item Degree $ks$ hypersurface in $\mathbb{P}(1,1,\cdots,1,s)$ satisfying assumptions in Theorem~\ref{thm_second_main_result} $(3)$. 
\end{enumerate}
\end{thm}
\subsection{Further questions}
Let $X$ be smooth hypersurface of degree $d$ in $\mathbb{P}(a_{0},a_{1},\cdots,a_{n})$ defined by $f$. Let $m=\mathrm{gcd}(d,\sum^{i=n}_{i=0}a_{i})$. We always reconstruct $$\bigoplus_{t\geq0} \Hom(\Delta,\Delta(tm)).\footnote{See Section 3 for the definition}$$
by the Hochschild--Serre algebra of $\Ku(X)$. For the examples in the paper, we can further reconstruct the sub-algebra $\bigoplus_{t\geq0} \mathrm{Jac}(f)_{tm}$ of the Jacobian ring $\mathrm{Jac}(f)$(If $X$ is the Veronese double cone, then $m=2$), and then reconstruct $X$. We ask the following natural questions,
\begin{Ques}
Let $X$ be smooth hypersurface of degree $d$ in $\mathbb{P}(a_{0},a_{1},\cdots,a_{n})$ defined by $f$.
 If $m=\mathrm{gcd}(d,\sum^{i=n}_{i=0}a_{i})\geq 2$, can sub-algebra $\bigoplus_{t\geq0} \mathrm{Jac}(f)_{tm}$ reconstruct $\mathrm{Jac}(f)$?   
\end{Ques}

In \cite{Donagi1986} and \cite{voisin2022schiffer}, the authors prove classical generic Torelli theorem for smooth hypersurface in $\mathbb{P}^n$ except for only finitely many cases. In our article, we prove the categorical Torelli theorem for smooth hypersurface of degree $d$, which does not divide $n+1$ and if $(d,n+1)=(4,4k+2)$ with $k\geq 150$. 
One of the remaining case is when $d$ divides $n+1$.
\begin{Ques}
Is it possible to reconstruct the Schiffer variation from the Kuznetsov component of degree $d$ hypersurface in $\mathbb{P}^n$ such that $d$ divides $n+1$?
\end{Ques}

\subsection{Related work}
Theorem \ref{thm_rennemo2023} is the main result of the unfinished preprint \cite{Rennemo2023hochschild}, which is referred to in some other papers, and which this paper replaces.
In \cite{huybrechts2016hochschild}, the authors prove a refined categorical Torelli theorem for cubic fourfold, whose Kuznetsov component is a $\mathrm{CY}2$ category, while we prove a refined categorical Torelli theorem for a Fano variety whose Kuznetsov component is a $\mathrm{CY}(2n+1)$ category. In \cite{lin2024serre}, the authors use the Hochschild--Serre algebra to prove the categorical Torelli theorem for a series of hypersurfaces and weighted hypersurfaces. In \cite{pirozhkov2022categorical}, the author proves the categorical Torelli theorem for a large class of hypersurfaces under the assumption that the equivalence of Kuznetsov components commutes with their respective rotation functors. In \cite{Dell2023categorical}, the authors prove the categorical Torelli theorem for the Veronese double cone under the same assumption above. In an upcoming preprint \cite{Lahoz2023categorical}, the authors prove categorical Torelli theorems following the line of \cite{pirozhkov2022categorical} for weighted hypersurfaces. 

\subsection{Organization of the paper}
In Section~\ref{section_intro_matrix_factorizations}, we recall the basics of the dg category of graded matrix factorizations associated with a weighted hypersurface $X$ of degree $d$ and dimension $n-1$; in particular, we recall that we can identify the dg enhancement $\Ku_{dg}(X)$ of the Kuznetsov component with a subcategory $\mathrm{Inj}_{\mathrm{coh}}(\mathbb{A}^{n+1},\mathbb{C}^*,\oh(d),\omega)$ of category of graded matrix factorization in localizing $\mathrm{Hqe}$ of dg categories. In Section~\ref{section_Serre_algebra}, we define the Hochschild--Serre algebra for an arbitrary smooth and proper dg category
and focus on the dg category of matrix factorization on affine LG model. In Section~\ref{section_categorical_Torelli_veronese_double_cone}, we prove Theorem~\ref{thm_main_result_categorical_Torelli_theorem}. In Section~\ref{section_IVHS_weighted_hypersurfaces}, we prove Theorem~\ref{thm_second_main_result}, Corollary~\ref{thm_the_third_main_result} and Theorem~\ref{theorem_refined_categorical_Torelli}.
In Section~\ref{section_improvement}, we prove Theorem~\ref{thm_actual_CTT}.

\subsection{Acknowledgements}
This material is based upon work supported by the National Science Foundation under Grant No. DMS-1928930 and by the Alfred P. Sloan Foundation under grant G-2021-16778, while JVR and SZ were in residence at the Simons Laufer Mathematical Sciences Institute (formerly MSRI) in Berkeley, California, during the Spring 2024 semester.
SZ is supported by ANR project FanoHK, grant ANR-20-CE40-0023, Deutsche Forschungsgemeinschaft under Germany's Excellence Strategy-EXC-2047$/$1-390685813.
Part of this research was performed while XL and SZ visited the Max Planck Institute for Mathematics and the Hausdorff Institute for Mathematics during the JTP program on Hodge theory, Cycles, and derived categories and MCM, Academy of Mathematics and System Sciences, China Academy of Science. They are grateful for excellent working conditions and hospitality. 
JVR is supported by the Research Council of Norway grant no.~302277.

\section{dg category of graded matrix factorizations}\label{section_intro_matrix_factorizations}
In this section, we recall the dg category of matrix factorizations, following \cite{BFK} and \cite{positselski_two_2011, efimov_coherent_2015}. 
We refer the reader to \cite{keller2006differential} for the basics of dg categories. 

Let $(X,G,L,\omega)$ be a quadruple where $X$ is a quasi-projective variety, $G$ is a reductive algebraic group acting on $X$, $L$ is a $G$-equivariant line bundle on $X$, and $\omega$ is a $G$-invariant section of $L$. 
Our main example is $(\mathbb{A}^{n+1},\mathbb{C}^{\ast},\mathcal{O}(d),\omega)$ with the action of $\mathbb{C}^*$ on $\mathbb{A}^{n+1}$ given by $\lambda\cdot(x_{0},x_{1},\cdots, x_{n})=(\lambda^{a_0} x_{0},\lambda^{a_1} x_{1},\cdots, \lambda^{a_n} x_{n})$; we will write $(\mathbb{A}^{n+1},\mathbb{C}^{\ast},\mathcal{O}(d),\omega)_{[a_0,a_1,\ldots,a_n]}$ to emphasise the weights when necessary. 
The line bundle $\cO(d)$ is the $\CC^*$-equivariant line bundle on $\AA^{n+1}$ whose $\CC^*$-invariant global sections are homogeneous degree $d$ polynomials, and so $\omega$ is given by such a polynomial.
We always assume that the $a_i$ are positive, and that $\omega$ only has an isolated singularity at $0\in \mathbb{A}^{n+1}$.

There is a dg category $\operatorname{Fact}(X,G,L,\omega)$, whose objects are quadruples $(\mathcal{E}_{-1},\mathcal{E}_{0},\Phi_{-1},\Phi_{0})$, where $\mathcal{E}_{-1}$ and $\mathcal{E}_{0}$ are $G$-equivariant quasi-coherent sheaves on $X$, and $\Phi_{-1}:\mathcal{E}_{0}\rightarrow \mathcal{E}_{-1}\otimes L$ and $\Phi_{0}:\mathcal{E}_{-1}\rightarrow \mathcal{E}_{0}$ are homomorphism of $G$-equivariant sheaves such that
\begin{align*}
    \Phi_{-1}\circ \Phi_{0}=\id_{\cE_0} \otimes \omega.\\
    (\Phi_{0}\otimes \id_L)\circ \Phi_{-1}=\id_{\cE_{-1}} \otimes \omega.
\end{align*}
The complex of morphisms between objects of $\operatorname{Fact}(X,G,L,\omega)$ is defined in \cite[Def.~3.1]{BFK}.
There is a full dg subcategory $\operatorname{Acyc}(X, G, L,\omega) \subseteq \Fact(X, G, L,\omega)$, analogous to the category of acyclic complexes in the category of complexes of sheaves. 

Both $\Fact$ and $\Acyc$ are pretriangulated dg categories, so their homotopy categories are naturally triangulated.
The absolute derived category is the Verdier quotient
\[
D^{\abs}(\operatorname{Fact}(X,G,L,\omega)) = [\operatorname{Fact}(X,G,L,\omega)]/[\operatorname{Acyc}(X,G,L,\omega)].
\]

Let $\operatorname{Inj}(X,G,L,\omega)\subseteq \operatorname{Fact}(X,G,L,\omega)$ be the dg sub-category whose objects are those where $\cE_{-1}$ and $\cE_0$ are injective $G$-equivariant quasi-coherent sheaves. 
The natural functor $[\operatorname{Inj}(X,G,L,\omega)]\rightarrow D^{\abs}(\operatorname{Fact}(X,G,L,\omega))$ is an equivalence of triangulated categories (see \cite[Proposition 3.11]{BFK} 
, and so $\Inj(X,G,L,\omega)$ is a dg enhancement of $D^{\abs}(\Fact(X,G,L,\omega))$.
   
Let $\operatorname{Inj_{coh}}(X, G, L,\omega)\subset \operatorname{Inj}(X, G, L,\omega)$ be the dg sub-category whose objects are quasi-isomorphic (that is, isomorphic in $D^{\abs}(\Fact(X, G, L,\omega))$) to objects with coherent components $\cE_i$.
  
In $\Fact(X,G,L,\omega)$ and all its subcategories, there is a shift auto-equivalence $[1]$, defined on objects by
\[
(\mathcal{E}_{-1},\mathcal{E}_{0},\Phi_{-1},\Phi_{0})[1] = (\mathcal{E}_{0},\mathcal{E}_{-1}\otimes L,-\Phi_{0},-\Phi_{-1}\otimes L).
\]
In the special case where $(X,G,L,\omega) = (\AA^{n+1},\CC^*,\cO(d), \omega)$, there is also a twist auto-equivalence
\[
\{1\}:\Fact(\mathbb{A}^{n+1},\mathbb{C}^{\ast},\mathcal{O}(d),\omega)\rightarrow \Fact(\mathbb{A}^{n+1},\mathbb{C}^{\ast},\mathcal{O}(d),\omega),
\]
which maps 
$$\xymatrix{\mathcal{E}_{-1}
\ar[r]^{\Phi_{0}}&\mathcal{E}_{0}\ar[r]^{\Phi_{-1}}&\mathcal{E}_{-1}(d)}$$
to 
$$\xymatrix{\mathcal{E}_{-1}(1)
\ar[r]^{\Phi_{0}(1)}&\mathcal{E}_{0}(1)\ar[r]^{\Phi_{-1}(1)}&\mathcal{E}_{-1}(d+1)}$$
We write $\{k\} = \{1\}^{\circ k}$ and $[k] = [1]^{\circ k}$.
It is then easy to see that we have an equivalence of functors $\{d\} \cong [2]$.

We consider the weighted projective space $\PP(a_0,\dots, a_n)$ as a Deligne--Mumford stack.
The function $\omega$ defines a smooth hypersurface substack $X \subseteq \PP(a_0,\dots, a_n)$.
 Let 
  $$\Ku(X):=\Big\langle \mathcal{O}_{X},\mathcal{O}_{X}(1),\cdots, \mathcal{O}_{X}\left(\sum^{n}_{i=0}a_{i}-d-1\right)\Big\rangle^{\perp}.$$
 Consider the natural enhancement $\operatorname{Inj_{coh}}(X)$, and let $\Ku_{dg}(X) \subset \operatorname{Inj_{coh}}(X)$ be the dg subcategory that enhances $\Ku(X)$. 
 Write $(1)$ for the quasi-endofunctor (of Fourier--Mukai type) of $\Ku_{dg}(X)$ that induces the degree shifting functor $(1):\Ku(X)\rightarrow \Ku(X)$ in the sense of \cite{huybrechts2016hochschild}. 

Let $\operatorname{Hqe(dg-cat)}$ be the localization of $\operatorname{dg-cat}$ with respect to quasi-equivalences of dg categories. 
The following theorem is \cite[Theorem 6.13]{BFK}, see also \cite{orlov2009derived, caldararu2679infinity}.
 
\begin{thm}\label{prop_Orlov_identification}\cite[Theorem 6.13]{BFK}
  There is an equivalence in $\operatorname{Hqe(dg-cat)}$,
  $$\Phi:\operatorname{\mathrm{Inj_{coh}}}(\mathbb{A}^{n+1},\mathbb{C}^{\ast},\mathcal{O}(d), \omega)_{[a_{0},a_{1},\cdots,a_{n}]}\cong \Ku_{dg}(X),$$
\end{thm}

\begin{prop}\cite[Thm.~1.2]{FeveroKellyLGfractionalcy}\label{serrefunctor}
The Serre functor of $[\operatorname{\mathrm{Inj_{coh}}}(\mathbb{A}^{n+1},\mathbb{C}^{\ast},\mathcal{O}(d), \omega)_{[a_{0},a_{1},\cdots,a_{n}]}]\simeq \Ku(X)$ is $( - \sum_{i=0}^n a_i)[n+1]$.
\end{prop}

\section{Hochschild--Serre algebra}\label{section_Serre_algebra}
\begin{thm}\label{functorbimodule}\cite{toen2007homotopy} Let $\mathcal{A}$ be a dg category over the field $k$. In $\operatorname{Hqe(dg-cat)}$, we have an equivalence
$$\mathcal{D}_{dg}(\mathcal{A}^{op}\otimes \mathcal{A})\cong \mathcal{R}\Hom_{c}(\mathcal{D}_{dg}(\mathcal{A}),\mathcal{D}_{dg}(\mathcal{A})),$$
where $\mathcal{R}\Hom_{c}(\mathcal{D}_{dg}(\mathcal{A}),\mathcal{D}_{dg}(\mathcal{A}))$ is the dg category of quasi-functors preserving coproducts.
\end{thm}
Suppose $\mathcal{A}$ is a smooth proper dg category. By \cite[Proposition 1.45]{Tab15}, there is a smooth proper dg algebra $A$, Morita equivalent to $\mathcal{A}$. In that case, the bimodules for the Serre functor and inverse of the Serre functor are constructed explicitly in \cite{shklyarov2007serre}, where the author defines the Serre functor (inverse Serre functor) of the triangulated category $\operatorname{Perf}(A)$. 
From now on, we interpret 
those bimodules as quasi-functors by Theorem~\ref{functorbimodule}. 

\begin{defn} The Hochschild cohomology and homology of a smooth proper dg category $\mathcal{A}$ are defined as
    \[
    \HH^{m}(\mathcal{A})=\Hom(\Id,\Id[m]) \text{ and } \HH_{m}(\mathcal{A})=\Hom(\Id,S[m]),
    \]
respectively, where $S$ is the Serre functor of $\cA$. 
\end{defn}
The Hochschild cohomology is an algebra, and Hochschild homology is a graded module over this algebra.
We now define an algebra that contains Hochschild cohomology and Hochschild homology and encodes the algebra structure of Hochschild cohomology and the module structure of Hochschild homology over Hochschild cohomology.

\begin{defn}(Hochschild--Serre algebra)
   Let $\mathcal{A}$ be a smooth, proper dg category. 
   Define a bi-graded algebra 
   $$\mathrm{HS}(\cA)=\bigoplus_{m,n\in \mathbb{Z}} \Hom(\Id,S^{m}[n])$$
 The multiplication map 
 $$\xymatrix@C=2cm{\Hom(\Id,S^{m_1}[n_1])\times \Hom(\Id,S^{m_2}[n_2])\ar[r]^{\times}&\Hom(\Id,S^{m_1+m_2}[n_1+n_2])}$$
 is defined as follows.
 For an element $(a,b)\in \Hom(\Id,S^{m_1}[n_1])\times \Hom(\Id,S^{m_2}[n_2])$, $a\times b$ is defined as the composition
 $$\xymatrix@C=3.5cm{\Id\ar[r]^{b}&\Id\circ S^{m_2}[n_2]\ar[r]^{a\circ \Id}&S^{m_1}[n_1]\circ S^{m_2}[n_2]=S^{m_1+m_2}[n_1+n_2]
 }.$$\end{defn}
It is shown in \cite[Definition 3.3]{lin2024serre} that $\mathrm{HS}(\cA)$ is an associative algebra.

Let $\operatorname{Hmo(dg-cat)}$ be the localization of $\operatorname{dg-cat}$ with respect to the Morita equivalences of dg categories. If $\mathcal{A}$ and $\mathcal{B}$ are smooth and proper, 
$\Hom_{\operatorname{Hmo(dg-cat)}}(\mathcal{A},\mathcal{B})$ is the set of isomorphism classes of objects in $\operatorname{Perf}(\mathcal{A}^{op}\otimes \mathcal{B})$ \cite[Corollary 1.44]{Tab15}, and the composition corresponds to tensor product. 
The following theorem is proved in \cite[Theorem 3.4]{lin2024serre}.

\begin{thm}\label{Moritaserrealgebra}
\cite[Theorem 3.4]{lin2024serre}\cite[Appendix A.4]{belmans2023hochschild}
If $\mathcal{A}\simeq\mathcal{B}$ in $\operatorname{Hmo(dg-cat)}$, then $\mathrm{HS}(\cA)\cong\mathrm{HS}(\cB)$.    
\end{thm}

\subsection{Examples of Hochschild--Serre Algebras}
This section gives examples of Hochschild--Serre algebras for many categories. 

\subsubsection{Orlov's algebra $\mathrm{HA}(X)$}\label{section_orlov_algebra}
Let $\cA=D^b(X)$ be the bounded derived category of coherent sheaves on a smooth projective variety $X$. In this case, $S_{\cA}=-\otimes\omega_X[l]$, where $l=\mathrm{dim}X$. Thus the Hochschild--Serre algebra $\mathrm{HS}(\cA)$ is given by 
\begin{align*}
\mathrm{HS}(\cA):=\bigoplus_{m,n\in\mathbb{Z}}\mathrm{Hom}(\mathrm{Id},\mathrm{S}_{D^b(X)}^m[n])
\cong&\bigoplus_{m,n\in\mathbb{Z}}\mathrm{Hom}_{D^b(X\times X)}(\iota_*\oh_X,\iota_*\omega_X^{\otimes m}[ml+n])\\
\cong&\bigoplus_{m,n\in\mathbb{Z}}\mathrm{Ext}_{X\times X}^{ml+n}(\iota_*\oh_X,\iota_*\omega_X^{\otimes m}),
\end{align*}
where $\iota:X\hookrightarrow X\times X$ be the diagonal inclusion. It is clear that $\mathrm{HS}(\cA)$ is isomorphic to the bi-graded algebra $\mathrm{HA}(X)$ in \cite{orlov2003derived}. 
For a given $m \in \ZZ$, the graded piece $\mathrm{HS}(\cA)(m,-ml)$ equals
\[
\mathrm{Hom}_{X\times X}(\iota_*\oh_X,\iota_*\omega_X^{\otimes m}) = \mathrm{Hom}_X(\oh_X,\omega_X^{\otimes m}),
\]
and summing over $m \ge 0$ exhibits the canonical ring of $X$ as a subring of $\mathrm{HS}(\cA)$.


\subsubsection{dg category of matrix factorization on affine LG model}
According to Theorem \ref{prop_Orlov_identification}, the Kuznetsov components of a hypersurface $X\subset\mathbb{P}(a_{0}, a_{1}, \cdots, a_{n})$, defined by a polynomial $\omega$ of degree $d$, is quasi-equivalent to the category of factorizations $\operatorname{Inj_{coh}}(\mathbb{A}^{n+1}, \mathcal{O}_{\mathbb{A}^{n+1}},\mathbb{C}^{\ast},\omega)_{[a_{0},\dots,a_{n}]}$. 
The Serre functor of $\operatorname{Inj_{coh}}(\mathbb{A}^{n+1},\mathcal{O}_{\mathbb{A}^{n+1}},\mathbb{C}^{\ast},\omega)_{[a_{0},a_{1},\dots,a_{n}]}$ is $-\otimes \mathcal{O}_{\mathbb{A}^{n+1}}(-\sum^{n}_{i=0}a_{i})[n+1]$ by Proposition \ref{serrefunctor}.

According to \cite{BFK}, the natural functors can be reinterpreted as Fourier--Mukai transformation of kernels, and the natural transformation between these functors is a morphism of kernels. 
We let $\Delta(m)$ be the kernel of the functor $-\otimes\mathcal{O}_{\mathbb{A}^{n+1}}(m)$. 
Let $\mu_{d}=\langle \exp{\frac{2\pi i}{d}}\rangle \subset \CC^*$, and note that $\mu_d$ acts on $\AA^{n+1}$ since $\CC^*$ does.
For any $g \in \mu_d$, let $W_{g}$ be the conormal bundle of the fix point locus $(\mathbb{A}^{n+1})^{g}$ in $\AA^{n+1}$, and let $k_{g} \in \ZZ$ be such that $\operatorname{det}(W_{g}) \cong \cO(k_g)$ as $\CC^*$-equivariant bundles.
Let $\omega_{g}:=\omega\vert_{(\mathbb{A}^{n+1})^{g}}$.
\begin{thm}\cite[Theorem 5.39]{BFK}\label{extendedhochschild}
 If $\omega$ has an isolated singularity exactly at $0\in \mathbb{A}^{n+1}$, then 
\begin{align*}
 \Hom(\Delta,\Delta(t)[m])\cong& \bigoplus_{\{g\in \mu_{d}\ \mid\  \operatorname{rk}W_{g} \equiv m\pmod2\}}\mathrm{Jac}(\omega_{g})_{t-k_{g}+d\left(\frac{m-\operatorname{rk}W_{g}}{2}\right)}
 \end{align*}
\end{thm}
\begin{rem}
We may unpack the quantities in the theorem above as follows.
Given $a_0,\dots, a_n, d$ and $g \in \mu_d$, we let $I_g = \{i \in \{0, \dots, n\} \mid g^{a_i} \not= 1\}$.
Then $\rk W_g = |I_g|$, $k_g = -\sum_{i \in I_g}a_i$, and $\omega_g$ is the image of $\omega$ in the ring $\CC[x_0,\dots, x_n]/(x_i)_{i \in I_g}$.
\end{rem}
\begin{rem}\label{compatibility}
In a similar way as with the Hochschild--Serre algebra, the direct sum
\[
\bigoplus_{t,m \in \ZZ} \Hom(\Delta, \Delta(t)[m])
\]
can be given a structure of a $\CC$-algebra where the multiplication is defined by the composition of functors.
Using the above theorem, we get a decomposition
\[
\Hom(\Delta, \Delta(t)[m]) = \bigoplus_{g \in \mu_d} \Hom(\Delta,\Delta(t)[m])_g.
\]
It is easy to check that $\bigoplus_{t \in \ZZ} \Hom(\Delta, \Delta(t))_1$ is closed under multiplication, and that the isomorphisms 
\[
\Hom(\Delta, \Delta(t))_1 \cong \Jac(\omega)_t,
\]
given by the theorem respect the multiplicative structures on both sides, i.e.~ we have ring homomorphisms 
\[
\Jac(\omega) \overset{\cong}{\to} \bigoplus_{t \in \ZZ} \Hom(\Delta, \Delta(t))_1 \into \bigoplus_{t,m \in \ZZ} \Hom(\Delta, \Delta(t)[m]).
\]
\end{rem}

\begin{lem}
\label{thm:recoveringSubJacobian}
Let weights $a_0, \dots, a_n > 0$ be given, and let $\omega, \omega'$ be two homogeneous degree $d$ polynomials, both with 0 as the unique critical point.
    Let $t = \gcd(d, \sum_{i=0}^n a_i)$,
     let
    \[
    \cA = \operatorname{Inj_{coh}}(\mathbb{A}^{n+1},\mathcal{O}_{\mathbb{A}^{n+1}},\mathbb{C}^{\ast},\omega)_{[a_{0},a_{1},\dots,a_{n}]} 
    \]
    and
    \[
    \cA' = \operatorname{Inj_{coh}}(\mathbb{A}^{n+1},\mathcal{O}_{\mathbb{A}^{n+1}},\mathbb{C}^{\ast},\omega')_{[a_{0},a_{1},\dots,a_{n}]}
    \]
   We assume the following:
    \begin{enumerate}
        \item The rings $\oplus_{i \ge 0}\Jac(\omega)_{it}$ and $\oplus_{i \ge 0}\Jac(\omega')_{it}$ are generated in degree $t$.
        \item For each $g \in \mu_d \setminus\{1\}$, we have either that $\rk W_g$ is odd, or
    \[
    \Jac(\omega_g)_{t-k_g-\frac{d\rk W_g}{2}} = 0.
    \]
    \end{enumerate}
    If 
    \[
    \cA \cong \cA',
    \]
    then we have an isomorphism of graded rings
    \[
    \bigoplus_{i \ge 0}\Jac(\omega)_{ti} \cong \bigoplus_{i \ge 0}\Jac(\omega')_{ti}.
    \]
\end{lem}
\begin{proof}
Using $S \cong \Delta(-\sum^{n}_{i=0}a_{i})[n+1]$ and $\Delta(d) \cong [2]$, we may find $k, l$ such that $\Delta(t) \cong S^k[l]$.
By our assumption and Theorem \ref{extendedhochschild}, we have that
\[
\Hom(\mathrm{Id}, S^k[l]) \cong \Hom(\Delta, \Delta(t)) = \Jac(\omega)_t.
\]
Using Remark \ref{compatibility}
, we find the subalgebra of $\mathrm{HS}(\cA)$ generated by $\Hom(\mathrm{Id}, S^k[l])$ is isomorphic to $\bigoplus_{i \ge 0} \Jac(\omega)_{it}$.
Since an equivalence $\cA \cong \cA'$ yields a bigraded isomorphism $\mathrm{HS}(\cA) \cong \mathrm{HS}(\cA')$ by Theorem \ref{Moritaserrealgebra}, the claim follows.
\end{proof}
Condition (1) holds if all $a_i$ are equal to 1.
Condition (2) is easily checked for given values of $(a_i)$ and $d$.
In the cases of interest to us, it can be checked by using the following lemma.
\begin{lem}
\label{thm:vanishingOfGPiece}
    Let $a_0,\dots, a_n$, $\omega$, $d$ and $t$ be as in Lemma \ref{thm:recoveringSubJacobian}.
    Assume that $a_i \le \frac{d}{2}$ for all $i$.
    Let $I = \{i \in \{0,\ldots, n\} \mid a_i = 1\}$
    If
    \[
    |I| > \frac{2t}{d-2},
    \]
    then $\mathrm{Jac}(\omega)_{t - k_g - \frac{d\rk W_g}{2}} = 0$ for all $g \in \mu_d \setminus \{1\}$.
\end{lem}
\begin{proof}
Let $g \in \mu_d \setminus \{1\}$, and let $J \subseteq \{0,\ldots, n\}$ be the set of $i$ such that $g^{a_i} \not= 1$.
Clearly, $I \subseteq J$.
We have
\[
k_g + d\rk W_g/2 = \sum_{i \in J} \left(\frac{d}{2} - a_i\right) \ge \sum_{i \in I} \left(\frac{d}{2} - a_i\right) = |I|\left(\frac{d}{2}-1\right),
\]
using here the assumption $a_i \le d/2$ for all $i$.
We then have 
\[
t - \left(k_g + \frac{d\rk W_g}{2}\right) \le t - |I|\left(\frac{d}{2} - 1\right) < 0,
\]
and the claim follows.
\end{proof}

\subsection{Smooth hypersurfaces in projective space}
We consider in this section the case where $a_0 = \dots = a_n =1$, so that $\omega$ is the defining equation of a degree $d$ hypersurface $Y$ of $\PP^n$.
We let $t = \gcd(d, n+1)$. In this case, the Kuznetsov component 
 $$\Ku(Y):=\Big\langle \mathcal{O}_{Y},\mathcal{O}_{Y}(1),\cdots, \mathcal{O}_{Y}\left(n-d)\right)\Big\rangle^{\perp},$$
 which is equivalent to $\mathrm{Inj}_{\mathrm{coh}}(\mathbb{A}^{n+1},\mathbb{C}^{\ast},\omega)_{[1,1,\dots,1]}$ by Theorem~\ref{prop_Orlov_identification}. 

 The following proposition along with the isomorphism $\Hom(\Delta, \Delta(i)[N]) \cong \Hom(\Delta, \Delta(i-d)[N+2])$ suffices to determine $\Hom(\Delta, \Delta(i)[j]$ for all $i,j$.
 

 \begin{prop}\label{HShypersurface}
Consider $\mathrm{Inj}_{\mathrm{coh}}(\mathbb{A}^{n+1},\mathcal{O}_{\mathbb{A}^{n+1}}, \mathbb{C}^{\ast},\omega)_{[1,1,\cdots,1]}$. Let $\Delta(i)$ be the kernel of functor $(i)$,

 \leavevmode\begin{enumerate}
 \item If $n+1$ is even, then
     \[
         \Hom(\Delta,\Delta(i)[0]) \cong \begin{cases}
         \mathrm{Jac}(\omega)_{i}, & i\neq \frac{(d-2)(n+1)}{2} \\
         \mathrm{Jac}(\omega)_{i}\oplus \mathbb{C}^{d-1},& i=\frac{(d-2)(n+1)}{2}
         \end{cases}
     \]
     and
     \[
     \Hom(\Delta,\Delta(i)[1]) = 0.
     \]
     \item If $n+1$ is odd, then  
     $$ \Hom(\Delta,\Delta(i)[0])\cong \mathrm{Jac}(\omega)_{i},$$
 and
 \[
 \Hom(\Delta,\Delta(i)[1]) \cong \begin{cases}
         \mathbb{C}^{d-1}, & i= \frac{(d-2)(n+1)}{2} - \frac{d}{2} \\
         0, & i\neq \frac{(d-2)(n+1)}{2} - \frac{d}{2}.
       \end{cases}
 \]
 \end{enumerate}
 \end{prop}

\begin{prop}
\label{thm:ordinaryHypersurfacesJacobian}
Assume that $3 \le d \le n$, and $(d,n+1) \not= (3,6)$. 
If $\omega$ and $\omega'$ are degree $d$ polynomials so that
\[
\operatorname{Inj_{coh}}(\mathbb{A}^{n+1},\mathcal{O}_{\mathbb{A}^{n+1}},\mathbb{C}^{\ast},\omega)_{[1,\dots,1]} \cong \operatorname{Inj_{coh}}(\mathbb{A}^{n+1},\mathcal{O}_{\mathbb{A}^{n+1}},\mathbb{C}^{\ast},\omega')_{[1,\dots,1]},
\]
then
\[
\oplus_{i\ge 0}\Jac(\omega)_{it} \cong \oplus_{i \ge 0}\Jac(\omega')_{it}.
\]
\end{prop}
\begin{proof}
First, suppose $d$ does not divide $n+1$. 
We then have $t \le d/2$.
Since then
\[
\frac{2t}{d-2} \le 1 + \frac{2}{d-2} \le 3 < n+1,
\]
the assumptions of Lemma \ref{thm:vanishingOfGPiece} are verified.
The claim then follows from Lemma \ref{thm:recoveringSubJacobian}.

Next, suppose $d$ divides $n+1$, we have  $\mathrm{gcd}(d,n+1)=d$, and $d\neq \frac{(d-2)(n+1)}{2}$ because 
the equality $d=\frac{(d-2)(n+1)}{2}$ holds for integers $3\leq d\leq n$ if and only if $(d,n+1)=(3,6)$. 
The claim follows from  Proposition~\ref{HShypersurface} and Lemma~\ref{thm:recoveringSubJacobian}. 
\end{proof}


Combining the above computation with the techniques of \cite{Don,CM_1990__73_2_121_0} gives the following consequence.
\begin{thm}\label{thm_rennemo2023}
 Let $X$ and $X'$ be smooth hypersurfaces of degree $d \le n$ in $\mathbb{P}^{n}$.
 Assume that $d$ does not divide $n+1$, and that $(d, n+1)$ is not of the form $(4,4k+2)$ for $k \in \ZZ$. 
 If $\Ku(X) \cong \Ku(X')$, then $X\cong X'$.
\end{thm}
\begin{proof}
The cases of $d \le 2$ or $n \le 2$ are trivial, so we may assume $d, n \ge 3$.
Let $\omega$ and $\omega'$ be the polynomials defining $X$ and $X'$, respectively, and let $t = \gcd(d,n+1)$.

We claim that there is up to $\GL(n+1)$-equivalence a unique surjective ring homomorphism $\oplus_{i \ge 0} S^{it}(\CC^{n+1}) \to \oplus_{i \ge 0} \Jac(\omega)_{it}$, and likewise for $\oplus_{i \ge 0} \Jac(\omega')$.

If $2t < d$, then this follows from \cite[Lem.~4.2]{Don}.
If $2t = d$, then since $(d,n+1)$ is not of the form $(4,4k+2)$, we must have $d > 4$.
Then \cite{CM_1990__73_2_121_0} gives the claim.

Since $\oplus_{i \ge 0} \Jac(\omega)_{it} \cong \oplus_{i \ge 0} \Jac(\omega')_{it}$ by Prop.~\ref{thm:ordinaryHypersurfacesJacobian}, we have that the degree $d$ parts of the Jacobian ideals of $\omega$ and $\omega'$ are $\mathrm{GL}(n+1)$-equivalent.
Thus by \cite[Prop.~1.1]{Don}, we must have $X \cong X'$.
\end{proof}

In the remaining case of $(d,n+1) = (4, 4k+2)$, we can apply Voisin's recent work \cite{voisin2022schiffer} to prove the following statement.
\begin{cor}\label{cor_generic_categorical_Torelli_quartic_hypesurface}
Let $Y$ and $Y'$ be generic smooth quartic hypersurfaces in $\mathbb{P}^{4k+1}$, with $k \ge 150$.
If there is a Fourier--Mukai type equivalence $\Phi:\Ku(Y)\simeq\Ku(Y')$, then $Y\cong Y'$. 
\end{cor}

\begin{proof}
By Proposition~\ref{thm:ordinaryHypersurfacesJacobian}, 
\[
\bigoplus_{i \ge 0}\mathrm{Jac}(Y)_{2i}\cong \bigoplus_{i \ge 0}\mathrm{Jac}(Y')_{2i}.
\]
By \cite[Section 2]{voisin2022schiffer}, we obtain $Y\cong Y'$. 
\end{proof}

\begin{rem}\label{rem_double_cover}
Consider a double cover $X\rightarrow\mathbb{P}^{4n+1}$ with branch divisor being the quartic hypersurface $Y\subset\mathbb{P}^{4n+1}$. By \cite[Corollary 4.6]{kuznetsov2019calabi}, the semi-orthogonal decomposition of $X$ is given by 
$$D^b(X)=\langle\Ku(X),\oh_X,\ldots,\oh_X(4n-1)\rangle,$$
and $\Ku(X)$ is a Calabi--Yau category such that $\mathcal{S}_{\Ku(X)}\cong [2n+1]$. The Rotation functor for $\Ku(X)$ is defined as $\bR:=\bL_{\oh_X}(-\otimes\oh_X(1))$. Then by \cite[Proposition 3.4, Corollary 3.18]{kuznetsov2019calabi}, $\bR^2\cong\tau[1]$, where $\tau$ is the geometric involution of $X$ coming from double cover. 
\end{rem}

\begin{cor}\label{cor_categorical_Torelli_double_cover}
Let $X,X'\rightarrow\mathbb{P}^{4n+1}$ be the double covers in Remark~\ref{rem_double_cover} such that $n\geq 150$.
Assume that there is a Fourier--Mukai type equivalence $\Phi:\Ku(X)\simeq\Ku(X')$ such that $\Phi$ commutes with rotation functors $\bR$ and $\bR'$ respectively. Then $X\cong X'$. 
\end{cor}

\begin{proof}
First note that by \cite[Proposition 7.10]{kuznetsov2017derived}, the equivariant category $\Ku_{\mathbb{Z}_2}(X)\simeq\Ku(Y)$ under the action of $\tau$. On the other hand, by Remark~\ref{rem_double_cover}, $\bR^2\cong\tau[1]$. Then $\Phi$ commutes with $\tau$ and $\tau'$ respectively. Then there is an induced Fourier--Mukai equivalence $\Phi^{\mathbb{Z}_2}:\Ku(Y)\simeq\Ku(Y')$, then by Corollary~\ref{cor_generic_categorical_Torelli_quartic_hypesurface}, we get $Y\cong Y'$ and in addition, the polynomials defining $Y$ and $Y'$ are projective equivalent, which in turn means that $X\cong X'$. 
\end{proof}

\begin{rem}\label{rem_refined_categorical_Torelli_series}
One could cook up a series of weighted hypersurfaces where the refined categorical Torelli theorem holds. Indeed, let $Y\subset\mathbb{P}^n$ be a degree $2d<n+1$ hypersurface such that $2d\not| n+1$. Then let $X\rightarrow\mathbb{P}^n$ be the double cover branched over the hypersurface $Y\subset\mathbb{P}^n$. Denote by $\tau$ the geometric involution on $X$ induced by the double cover. Then by \cite[Proposition 3.4, Corollary 3.18]{kuznetsov2019calabi}, $\bR^d[-1]\cong\tau$, where $\bR$ is the rotation functor of $\Ku(X)$. Then refined categorical Torelli theorem holds for $X$ by a similar argument in Corollary~\ref{cor_categorical_Torelli_double_cover}, from Theorem~\ref{thm_rennemo2023}.
\end{rem}

\section{Categorical Torelli theorem for Veronese double cone}\label{section_categorical_Torelli_veronese_double_cone}
In this section, we prove Theorem ~\ref{thm_main_result_categorical_Torelli_theorem}. 
First, we show the Hochschild--Serre algebra of Kuznetsov component $\Ku_{dg}(Y_1)$ of the Veronese double cone $Y_1$ reconstructs the even degree part of its Jacobian ring. 
Then, we prove that the categorical Torelli theorem holds for any Veronese double cone by extending the isomorphism of the even degree component of the Jacobian ring to all degrees. 

We start by showing that the defining equation of $Y_1$ can be written as 
$\omega=x^{2}_{5}+\omega^{b}$, where $\omega^{b}=x^{3}_{4}+\Phi_{1}x_{4}+\Phi_{2},$
and $\Phi_1,\Phi_2$ are degree $4$ and degree $6$ homogeneous polynomials respectively. 

\begin{lem}
Let $Y^{b}_{1}\subset \mathbb{P}(1,1,1,2)$ be a smooth hypersurface of degree $6$. 
After a variable change, the equation of $Y^{b}_{1}$ can be assumed to be of the form $\omega^{b}=x^{3}_{4}+\Phi_{1}(x_{1},x_{2},x_{3})x_{4}+\Phi_{2}(x_{1},x_{2},x_{3})$, where $\Phi_{1}$ is a homogeneous degree $4$ polynomial, and $\Phi_{2}$ is a homogeneous degree $6$ polynomial. 
\end{lem}
\begin{proof}
  We write $\omega^{b}=\alpha\cdot x^{3}_{4}+\Phi(x_{1},x_{2},x_{3})x^{2}_{4}+\Phi'_{1}(x_{1},x_{2},x_{3})x_{4}+\Phi'_{2}(x_{1},x_{2},x_{3})$.
  \begin{itemize}
        \item If
        $\alpha=0$. Then $\omega^{b}=\Phi x^{2}_{4}+\Phi'_{1}x_{4}+\Phi'_{2}$. Take  $[0,0,0,1]\in Z(\omega^{b})\subset \mathbb{C}^{4}$, it is also a common solution of $\frac{\partial \omega^{b}}{\partial x_{4}}=2x_{4}\Phi+\Phi'_{1}$, and $\frac{\partial \omega^{b}}{\partial x_{i}}=\frac{\partial \Phi}{\partial x_{i}}x^{2}_{4}+\frac{\Phi'_{1}}{\partial x_{i}}x_{4}+\frac{\partial \Phi'_{2}}{\partial x_{i}}, i=1,2,3$. Then $Y_{1}$ is not smooth. 
      \item Therefore $\alpha\neq 0$. We assume $\alpha=1$, take a coordinate change (automorphism of $\mathbb{P}(1,1,1,2)$) $x_{i}\mapsto x_{i},\ i=1,2,3;\ x_{4}\mapsto x_{4}+\Phi'(x_{1},x_{2},x_{3})$, where $2\Phi'(x_{1},x_{2},x_{3})=-\Phi(x_{1},x_{2},x_{3})$, then the equation of $\omega^{b}$ can be written as the form $x^{3}_{4}+\Phi_{1}(x_{1},x_{2},x_{3})x_{4}+\Phi_{2}(x_{1},x_{2},x_{3})$.
    
  \end{itemize}
  Thus, the equation of $Y^{b}_{1}$ is of the form $\omega^{b}=x^{3}_{4}+\Phi_{1}(x_{1},x_{2},x_{3})x_{4}+\Phi_{2}(x_{1},x_{2},x_{3})$.
  \end{proof}
\begin{lem}\label{varoneseform}
  Let  $Y_{1}$ be the double cover of $\mathbb{P}(1,1,1,2)$ branched along a smooth degree $6$ hypersurface. 
  After a variable change, the equation defining $Y_{1}$ can be assumed to be of the form $\omega=x^{2}_{5}+x^{3}_{4}+\Phi_{1}(x_{1},x_{2},x_{3})x_{4}+\Phi_{2}(x_{1},x_{2},x_{3})$. Furthermore, $\mathrm{Jac}(\omega)=\mathrm{Jac}(\omega^{b})$.   
\end{lem}
\begin{proof}
The equation of $Y_{1}$ may have a term $\Phi(x_{1},x_{2},x_{3},x_{4})x_{5}$. Take a coordinate change $x_{i}\mapsto x_{i},\ i=1,2,3,4;\ x_{5}\mapsto x_{5}+\Phi'(x_{1},x_{2},x_{3},x_{4})$, where $2\Phi'(x_{1},x_{2},x_{3},x_{4})=-\Phi(x_{1},x_{2},x_{3},x_{4})$. Then the equation defining $Y_{1}$ is of the form $\omega=x^{2}_{5}+x^{3}_{4}+\Phi_{1}(x_{1},x_{2},x_{3})x_{4}+\Phi_{2}(x_{1},x_{2},x_{3})$. Furthermore, 
$$ \mathrm{Jac}(\omega)\cong \frac{\mathbb{C}[x_{1},x_{2},x_{3},x_{4},x_{5}]}{\langle2x_{5},\frac{\partial\omega^{b}}{\partial x_{1}}, \frac{\partial\omega^{b}}{\partial x_{2}}, \frac{\partial\omega^{b}}{\partial x_{3}}, \frac{\partial\omega^{b}}{\partial x_{4}}\rangle}
             \cong \frac{\mathbb{C}[x_{1},x_{2},x_{3},x_{4}]}{\langle\frac{\partial\omega^{b}}{\partial x_{1}}, \frac{\partial\omega^{b}}{\partial x_{2}}, \frac{\partial\omega^{b}}{\partial x_{3}}, \frac{\partial\omega^{b}}{\partial x_{4}}\rangle} 
             = \mathrm{Jac}(\omega^{b}).$$

\end{proof}

\begin{thm}\label{lemma_even_degree_piece_Jacobian_ring}
The Kuznetsov component $\Ku(Y)$ of the Veronese double cone $Y$ reconstructs the even degree components of its Jacobian ring $\mathrm{Jac}(\omega)$.
\end{thm}   
\begin{proof}
According to Lemma~\ref{varoneseform},
condition (1) of Lemma \ref{thm:recoveringSubJacobian} is easy to verify. Condition (2) is also held by simple calculation. Then, combine Lemmas \ref{thm:recoveringSubJacobian} and \ref{thm:vanishingOfGPiece}, and the results follow.

\end{proof}

\subsection{Categorical Torelli theorem for Veronese double cones}
We break the proof of Theorem~\ref{thm_main_result_categorical_Torelli_theorem} into a sequence of lemmas and propositions. In this section, all homomorphisms and isomorphisms will be of graded $\CC$-algebras.
We will write $=$ to denote a canonical isomorphism.

Let $\cC$ be the category of non-negatively graded $\CC$-algebras, and let $d \ge 0$.
There is a functor $\tau_d \colon \cC \to \cC$ given by $\tau_d(R) = R/R_{> d}$. It admits a left adjoint $\nu_d \colon \cC \to \cC$, defined as follows.
Let $R \in \cC$, and let $\bar R = \CC[x_i]_{i \in S}$, where $S \subset R$ is the set of homogeneous elements.
There is an obvious graded $\CC$-algebra homomorphism $\phi \colon \bar R \to R$, and we let $I \subset \bar R$ be the ideal generated by $\ker \phi \cap \bar R_{\le d}$.
Finally set $\nu_d(R) = \bar R/I$. 

\begin{rem}\leavevmode
	\label{thm:genRels}
 \begin{enumerate}
     \item If $R \in \cC$ is defined as a $\CC$-algebra by homogeneous generators $x_i$ and homogeneous relations $f_j$, then $\nu_d(R)$ is canonically isomorphic to the $\CC$-algebra generated by the $x_i$ of degree $\le d$, modulo the relations $f_j$ of degree $\le d$.
     \item Given $R \in \cC$, we write $R^e \subseteq R$ for the even-degree subalgebra.
There is a canonical homomorphism $\nu_d(R^e) \to \nu_d(R)^e$ which is not in general an isomorphism, e.g. take $d = 1$ and $R = \CC[x]$ with $x$ of degree 1.
\item Let $\phi \colon R \to S$ be a homomorphism of nonnegatively graded $\CC$-algebras. Then $\nu_d(\phi) \colon \nu_d(R) \to \nu_d(S)$ is an isomorphsm if and only $\phi$ maps $R_{\le d}$ bijectively to $S_{\le d}$.
 \end{enumerate}

\end{rem}






\subsection{Passing from even degrees to all degrees}
The point of this section is Proposition \ref{thm:liftingFromEven}.

Consider the graded polynomial ring $\CC[x_1,x_2,x_3,y]$, where the $x_i$ have degree 1 and $y$ has degree 2.
Let $F \in \CC[x_i, y]$ be a polynomial of the form
\[
F = y^3 + g(x_1,x_2,x_3)y + h(x_1,x_2,x_3),
\]
where $g$ is of degree 4 and $h$ is of degree 6.
We let $J$ by
\[
J = Jac(F) = \CC[x_i,y]/(\partial_y F, \partial_{x_1} F, \partial_{x_2} F, \partial_{x_3} F).
\]
Since $\partial_{x_i}F$ has degree $5$, by Remark~\ref{thm:genRels} (1) we get
\begin{equation}
	\label{eqn:genRelsNu4}
\nu_4(J) = \CC[x_1,x_2,x_3,y]/(3y^2 + g),
\end{equation}
\begin{lem}
Let $\phi \colon \CC[a_{ij}, y]_{1 \le i,j \le 3} \to \nu_4(J)^e$ be the homomorphism defined by $\phi(a_{ij}) = x_ix_j$ and $\phi(y) = y$.
Then $\phi$ is surjective with $\ker \phi$ generated by
\begin{gather*}
	a_{ij} - a_{ji}\ \ \ \ \ 1 \le i, j \le 3 \\
	a_{ij}a_{kl} - a_{il}a_{kj}\ \ \ \ \ 1 \le i,j,k,l \le 3 \\
	3y^2 + \bar g,
\end{gather*}
where $\bar g$ denotes a lift of $g$ along $\CC[a_{ij}] \to \CC[x_1,x_2,x_3]$.
\end{lem}
\begin{proof}
	It is well known that the analogous homomorphism $\CC[a_{ij}] \to \CC[x_1,x_2,x_3]^e$ is surjective with kernel generated by the first two sets of relations.
	The claim of the lemma follows easily from this.
\end{proof}
By the above Lemma, $\nu_4(J)^e$ is defined by generators and relations of degree $\le 4$, and so we have by Remark~\ref{thm:genRels} (1) a canonical isomorphism
\begin{equation*}
	\nu_4(J)^e = \nu_4(\nu_4(J)^e).
\end{equation*}
Since the homomorphism $\nu_4(J)^e \to J^e$ is an isomorphism in degrees $\le 4$, by Remark~\ref{thm:genRels} (3), we get an isomorphism
\begin{equation*}
	\nu_4(\nu_4(J)^e) = \nu_4(J^e),
\end{equation*}
and so finally we find
\begin{equation}
\nu_4(J)^e = \nu_4(J^e).
\end{equation}

\subsubsection{Geometric interpretation of $\nu_4(J)$}
Let $\PP = \PP(1,1,1,2)$.
We have a canonical isomorphism
\begin{equation}
	\label{eqn:coordinateRingP}
\CC[x_1,x_2,x_3,y] = \bigoplus_{n \ge 0} H^0(\PP,\cO(n)),
\end{equation}
Let now $S \subset \PP$ be the surface defined by $3y^2 + g$. Since $S$ doesn't pass through the unique singular point $[0, 0, 0, 1]$ of $\mathbb{P}$, the geometry of $S$ is similar to a hypersurface of projective space. For example, $\mathcal{O}_{S}(n)$ is invertible for any $n$, and the dualizing sheaf $\omega_{S}\cong \mathcal{O}_{S}(4-5)=\mathcal{O}_{S}(-1)$. We refer the reader to \cite{Mori1975}.
\begin{lem}
	We have 
	\[
	\nu_4(J) = \bigoplus_{n \ge 0} H^0(S, \cO(n)).
	\]
\end{lem}
\begin{proof}
This follows from \eqref{eqn:genRelsNu4} and the short exact sequences
\[
0 \to H^0(\PP, \cO_{\PP}(d-4)) \overset{\cdot 3y^2 + g}\to H^0(\PP, \cO_{\PP}(d)) \to H^0(S, \cO_S(d)) \to 0.
\]
\end{proof}

It follows that we have canonical isomorphisms
\begin{equation}
\nu_4(J^e) = \nu_4(J)^e = \bigoplus_{n \ge 0} H^0(S, \cO(2n)).
\end{equation}

\subsubsection{Comparing different $F$}
Let $F'$ be a polynomial of the same shape as $F$,
\[
F' = y^3 + g'y + h',
\]
and write $J'$, $R'$, $S'$ for the objects defined similarly to $J, R, S$.
\begin{lem}
\label{thm:extendGeometricEvenToFullAlgebra}
	Every isomorphism of graded $\CC$-algebras
	\[
	\bigoplus_{n \ge 0} H^0(S, \cO_S(2n)) \cong \bigoplus_{n \ge 0} H^0(S', \cO_{S'}(2n))
	\]
	is the restriction of some isomorphism of graded $\CC$-algebras
	\[
	\bigoplus_{n \ge 0} H^0(S, \cO_S(n)) \cong \bigoplus_{n \ge 0} H^0(S', \cO_{S'}(n)).
	\]
\end{lem}
\begin{proof}
	Since $S$ and $S'$ are obtained by taking the Proj of the first two rings, we get an induced isomorphism of schemes
	\[
	f \colon S \to S'
	\]
	together with an isomorphism of line bundles $\beta \colon f^*\cO_{S'}(2) \cong \cO_{S}(2)$.
	
	By adjunction, we have isomorphisms $\alpha \colon \omega_S \to \cO_S(-1)$ and $\alpha' \colon \omega_S' \to \cO_{S'}(-1)$.
	We also have a canonical isomorphism $df \colon f^*(\omega_S') \to \omega_S$.
	Rescaling $\alpha$ if necessary, we can (since $H^0(S,\cO_S) = \CC$) get an equality of isomorphisms
	\[
	\beta \circ f^*(\alpha'^{\otimes -2}) = \alpha^{\otimes -2} \circ df^{\otimes -2} \colon f^*(\omega_S')^{\otimes -2} \to \cO_S(2)
	\]
	The induced isomorphism
	\[
	\bigoplus_{n \ge 0} H^0(S, \cO_S(n)) \overset{\alpha^{-1}}{\to} \bigoplus_{n \ge 0} H^0(S, \omega_S^{\otimes -n}) \overset{df}{\to} \bigoplus_{n \ge 0} H^0(S', \omega_{S'}^{\otimes -n}) \overset{\alpha'}{\to} \bigoplus_{n \ge 0} H^0(S',\cO_{S'}(n))  
	\]
	then agrees with the given one upon restriction to $\oplus_{n \ge 0} H^0(S, \cO(2n))$.
\end{proof}

\begin{cor}
	\label{thm:nu4JevenLifting}
	Let $\psi \colon \nu_4 (J^e) \to \nu_4 (J'^e)$ be an isomorphism.
	Then there exists an isomorphism $\chi \colon \nu_4 (J) \to \nu_4 (J')$ which induces $\psi$ by restriction.
\end{cor}
\begin{proof}
	This follows from the Lemma \ref{thm:extendGeometricEvenToFullAlgebra} along with $\nu_4(J^e) = \nu_4(J)^e$, and $\nu_4(J) = \oplus H^0(S, \cO(n))$, along with the primed versions.
\end{proof}

\begin{prop}
	\label{thm:liftingFromEven}
	Let $\phi \colon J^e \to J'^e$ be an isomorphism of graded $\CC$-algebras.
	There is an isomorphism of graded $\CC$-algebras $\psi \colon \CC[x_i, y] \to \CC[x_i, y]$, whose restriction to $\CC[x_i,y]^e \to \CC[x_i,y]^e$ induces the isomorphism $\phi \colon J^e \to J'^e$.
\end{prop}
\begin{proof}
	We have the following maps of Hom sets:
	\begin{align*}
		A \colon \Hom(J^e, J'^e) &\to \Hom(\nu_4(J^e), \nu_4 (J'^e)) \\
		B \colon \Hom(\nu_4(J), \nu_4(J')) &\to \Hom(\nu_4(J^e), \nu_4 (J'^e)) \\
		C \colon \Hom(\nu_4J, \nu_4J') &\to \Hom(\nu_2 (J), \nu_2 (J')) = \Hom(\CC[x_i,y], \CC[x_i,y]),
	\end{align*}
where the last line uses $\nu_2(J) = \nu_2(J') = \CC[x_i,y]$.

By Corollary \ref{thm:nu4JevenLifting}, $B$ is surjective, so we can choose a  $\widetilde{\psi}$ such that $B(\widetilde{\psi}) = A(\phi)$.
We then let $\psi = C(\widetilde \psi)$.
The restriction of and tracing through the definitions of the maps involved shows that $\phi$ is induced by $\psi$ in the way claimed.
\end{proof}

\subsubsection{Projective equivalence}
\begin{lem}\label{Matheyauweighted}
	\label{thm:upgradedDonagi}
	Let $n_1,\ldots, n_k$ be non-negative integers, and let 
	\[
	R = \CC[x_i^{(d)}]_{\substack{1 \le d \le k \\ 1 \le i \le n_d}}
	\]
	where $x_i^{(d)}$ is homogeneous of degree $d$.
	Let $F, F' \in R_D$, and assume that we have an equality of vector subspaces of $R_D$:
	\[
	\left\langle x^{(d)}_i \frac{\partial}{\partial x^{(d)}_j} F \right\rangle_{\substack{1 \le d \le k \\ 1 \le i,j \le n_d}} = \left\langle x^{(d)}_i \frac{\partial}{\partial x^{(d)}_j} F' \right\rangle_{\substack{1 \le d \le k \\ 1 \le i,j \le n_d}}.
	\]
	Then there exists an automorphism $\phi \colon R \to R$ such that $\phi(F) = F'$.
\end{lem}
\begin{proof}
This is a straightforward generalization of the proof of \cite[Proposition 1.1]{Don}, which is the special case $k = 1$.
Take $G = \prod_{d=1}^k GL(n_d)$, which naturally acts on the ring $R$ by automorphisms, and so on the vector space $R_D$.
The subspaces in the statement of the lemma are identified with $T_F(GF)$ and $T_F'(GF')$, respectively.
The elements $F, F' \in R_D$ satisfy the hypothesis of \cite[Lemma 1.2.]{Don}, so they lie in the same $G$-orbit, and in particular are related by an automorphism of $R$.
\end{proof}

\begin{prop}\label{eventoall}
	Let $\phi \colon J^e \to J'^e$ be an isomorphism.
	There is an automorphism of $\CC[x_i, y]$ sending $F$ to $F'$.
\end{prop}
\begin{proof}
	By Proposition \ref{thm:liftingFromEven}, there exists an isomorphism $\psi \colon \CC[x_i,y] \to \CC[x_i,y]$ inducing $\phi$.
	Replacing $F'$ by $\psi^{-1}(F')$, we then have that
	\[
	(\partial_{x_1} F, \partial_{x_2} F, \partial_{x_3} F, \partial_{y} F) = (\partial_{x_1} F', \partial_{x_2} F', \partial_{x_3} F', \partial_{y} F').
	\]
	Comparing the degree 6 pieces of these ideals, we find that
	\begin{align*}
	&\langle x_ix_j(3y^2 + g), 3y^3 + yg, x_iy\partial_{x_j}g + x_i\partial_{x_j}h \rangle_{i,j} \\
	= &\langle x_ix_j(3y^2 + g'), 3y^3 + yg', x_iy\partial_{x_j}g' + x_i\partial_{x_j}h' \rangle_{i,j} \subseteq \CC[x_i,y]_6.
\end{align*}
Restricting to the subspaces of polynomials with no $y^2$-term, we get
\[
\langle 3y^3 + yf, x_iy\partial_{x_j}f + x_i\partial_{x_j}g \rangle_{i,j} = \langle 3y^3 + yf', x_iy\partial_{x_j}f' + x_i\partial_{x_j}g' \rangle_{i,j},
\]
or equivalently
\[
\langle y\frac{\partial}{\partial y} F, x_i\frac{\partial}{\partial x_j} F \rangle_{i,j} = \langle y\frac{\partial}{\partial y} F', x_i\frac{\partial}{\partial x_j} F' \rangle_{i,j}
\]
We then conclude by Lemma \ref{thm:upgradedDonagi}.
\end{proof}
\subsubsection{Categorical Torelli theorem for Veronese double cones}
\begin{thm}\label{thm_categorical_Torelli_theorem_for_Veronese_double_cone}
Let $Y$ and $Y'$ be the degree $6$ smooth hypersurface of $\mathbb{P}(1,1,1,2,3)$. If $\Ku(Y)\simeq\Ku(Y')$, then $Y\cong Y'$.     
\end{thm}
\begin{proof}
 By Theorem~\ref{lemma_even_degree_piece_Jacobian_ring}, we get an isomorphism of the even degree of graded Jacobian rings,
 $$\bigoplus_{i}\mathrm{Jac}(Y)_{2i}\cong \bigoplus_{i}\mathrm{Jac}(Y')_{2i}.$$
 Then by Proposition~\ref{eventoall}, the isomorphism of the even degree of graded Jacobian ring extends to the isomorphism of graded Jacobian rings, and we have an automorphism of $\mathbb{P}(1,1,1,2,3)$ inducing isomorphism $Y\cong Y'$. 
\end{proof}
\section{IVHS for (weighted) hypersurfaces and generic Torelli}\label{section_IVHS_weighted_hypersurfaces}
Let $X$ be a smooth projective variety over $\mathbb{C}$. It admits the Hodge decomposition $$H^n(X,\mathbb{C})\cong\bigoplus_{p+q=n}H^q(X,\Omega_X^p),$$
where $H^q(X,\Omega_X^p)$ is identified with space of harmonic $(p,q)$-forms. 
The Hodge structure of $V_{\mathbb{C}}:=H^*(X,\mathbb{C})$ consists of a Hodge filtration 
$$0\subset F^n(V_{\mathbb{C}})\subset F^{n-1}(V_{\mathbb{C}})\subset\ldots\subset F^0(V_{\mathbb{C}})\subset V_{\mathbb{C}},$$
such that there is a decomposition $V_{\mathbb{C}}=F^p(V_{\mathbb{C}})\oplus\overline{F^{n+1-p}(V_{\mathbb{C}})}$ for every $p$ with $0 \le p \le n$. The cup product of the cohomology induces a bilinear pairing, which leads to the notion of the polarized Hodge structure. Denote by $D$ the moduli space of such polarized Hodge structure of all weights and denote by $\mathrm{Def}(X)$ the space of deformation space of $X$ if it exists, which can be regarded as moduli space of complex structures on $X$. Let $\langle ,\rangle$ be the nondegenerate form of $H^{n}_{\mathbb{Z}}=H^{n}(X,\mathbb{Z})_{tf}$, and $G_{\mathbb{Z}}:= Aut(H^{n}_{\mathbb{Z}}, \langle , \rangle)$. Thus we get a map 
$$\mathcal{P}:\mathrm{Def}(X)\rightarrow\frac{D}{G_{\mathbb{Z}}},$$
called period map, which associates a polarized Hodge structure to a complex structure. The (generic) Torelli problem asks if $\mathcal{P}$ is generically injective.
\par
Taking the differential $d\mathcal{P}$ of the period map and the following theorem tells us what $d\mathcal{P}$ looks like.

\begin{prop}\cite{carlson2017period}
The derivative of the period map $\mathcal{P}$ induces a map

$$ d\mathcal{P} : H^1(X,T_X) \longrightarrow \bigoplus_{p+q=n}\mathrm{Hom}(H^p(X,\Omega_X^q),H^{p+1}(X,\Omega_X^{q-1})).$$
More precisely, given element $\xi\in H^1(X,T_X)$, it  defines a map from $H^p(X,\Omega_X^q)$ to $H^{p+1}(X,\Omega_X^{q-1})$ by contracting vector field. 
\end{prop}
Infinitesimal Variation of Hodge structures(IVHS) \cite{carlson2017period} consists of integral polarised Hodge structure, which is transcendental data and 
$$\{H^1(X,T_X),\bigoplus_{p+q=n}\mathrm{Hom}(H^p(X,\Omega_X^q),H^{p+1}(X,\Omega_X^{q-1})), d\mathcal{P}\}$$ as algebraic data. Let $H_{\mathbb{Z}}$ be a free abelian group of finite rank. Let 
$$H_{\mathbb{Z}}\otimes \mathbb{C}\cong \bigoplus^{p=n}_{p=0}H^{n-p,p}.$$
be a polarized weight $n$ Hodge structure of $H_{\mathbb{Z}}$.
We refer the reader to \cite{carlson2017period} for the abstract definition of $\mathrm{IVHS}$. From now on, we refer to the algebraic data if we talk about $\mathrm{IVHS}$ of Hodge structure of $H_{\mathbb{Z}}$.
\begin{defn}\footnote{In \cite[Section 3]{Saio1986}, the definition of $\mathrm{IVHS}$ consist of an extra data compared to the definition in \cite{Donagi1986}, namely the polarization $Q: H^{n-p,p}\times H^{p,n-p} \rightarrow \mathbb{C}$.}
 An $\mathrm{IVHS}$ of the Hodge structure of $H_{\mathbb{Z}}$ consist of 
 \begin{itemize}
     \item A collection of vector spaces over $\mathbb{C}$, $\{T, H^{n-p,p}\}_{p=0,1,\cdots,n}$.
     \item A collection of bilinear maps,
 $$\delta_{p}:T\times H^{n-p,p}\rightarrow H^{n-p-1,p+1},$$
 \end{itemize}
\end{defn}


\begin{defn}
    Given two $\mathrm{IVHS}$, $\{T, H^{n-p,p},\delta_{p}\}_{p=0,1,\cdots,n}$ and $\{T', H'^{n-p,p},\delta'_{p}\}_{p=0,1,\cdots,n}$, we say they are isomorphic if there are isomorphisms $T \cong T'$, $H^{n-p,p} \cong H'^{n-p,p}$ for all $p$, and commutative diagrams,
    $$\xymatrix{T\times H^{n-p,p}\ar[r]^{\delta_{p}}\ar[d]&H^{n-p-1,p+1}\ar[d]\\
    T'\times H'^{n-p,p}\ar[r]^{\delta'_{p}}&H'^{n-p-1,p+1}\\}$$
 \end{defn}

 Given a smooth, projective variety $X$ of dimension $n$, we obtain an IVHS by setting $T = H^1(X,T_X)$, $H^{n-p,p} = H_{prim}^{p}(X,\Omega_{X}^{n-p})$, with the multiplication map induced by contracting vector field.
 


The variational Torelli problem asks if the algebraic data of infinitesimal variation of Hodge structure determines $X$ up to isomorphism.

\subsection{Infinitesimal Variation of Hodge structures of weighted hypersurfaces}\label{IVHS_Generalized_double_cone_section}
Let $Y$ be a smooth hypersurface of the weighted projective space $\mathbb{P}(a_{0}, a_{1}, \cdots,  a_{n})$ defined by a degree $d$ homogeneous polynomial, where $gcd(a_{0}, a_{1}, \cdots, a_{n})=1$. We assume the dimension of $Y$ is greater than $2$. 
Let $\mathrm{Jac}(f)$ be the Jacobian ring defined by $f$, and let $w=\sum^{n}_{i=0} a_{i}$. The following lemmas show the $\mathrm{IVHS}$ of $Y$ are graded components of its Jacobian ring. 

\begin{lem}\label{Groresidul}
Let $Y$ be the (weighted) hypersurface as above.  Let $H^{\ast}_{prim}(Y)$ be the primitive cohomology group. Then there are isomorphisms 
   $$H^{p,n-p}_{prim}(Y)\cong\mathrm{Jac}(f)_{(n-p+1)d-w},$$ and
   $$H^{1}(Y,T_{Y})\cong\mathrm{Jac}(f)_{d}.$$
   which are compatible with multiplication.
\end{lem}
\begin{proof}
    The first isomorphism follows from \cite[Theorem 1.2]{Tu1986}. Since the dimension of $Y$ is greater than $2$, the second isomorphism follows from \cite[Proposition 2.3]{Saio1986}.
\end{proof}
\begin{lem}\label{weightedIVHS}
Let $Y$ be the (weighted) hypersurface as above. The $\mathrm{IVHS}$ of $Y$ is isomorphic to the one obtained by setting
$$T = \mathrm{Jac}(f)_{d},\ H^{p, n-p} = \Jac(f)_{(n-p+1)d-w},$$
with bilinear maps defined via multiplication in the Jacobian ring,
$$\mathrm{Jac}(f)_{d}\times \mathrm{Jac}(f)_{td-w}\rightarrow \mathrm{Jac}(f)_{(t+1)d-w},\ \ \ t=1,\ 2,\ \cdots,\ n.$$
\end{lem}
\begin{proof}
    This follows from Lemma \ref{Groresidul}.
\end{proof}
\begin{rem}\label{weightedcommondivisor}
  Let $m=\mathrm{gcd}(d,w)$, then $d,\ td-w$ are all multiple of $m$ for $t= 1,\ 2,\ \cdots,\ n+1$.  
\end{rem}
Let $Z\subset\mathbb{P}(\underbrace{1,\ldots,1}_{r+1},2n,kn)$ be the weighted hypersurface of degree $2kn$ such that $k\geq 3$ is odd and $r\geq 2$. We call $Z$ the generalized Veronese double cone.
\begin{prop}\cite[Theorem B]{Saio1986}\label{IVHSgeneric}
 If $m<2n$, the $\mathrm{IVHS}$ of $Z$ recovers $Z$ generically.
\end{prop}

Let $X$ be a smooth hypersurface of degree $ks$ in $\mathbb{P}(\underbrace{1,1,\cdots,1}_{r+1},s)$ defined by $\omega$. Assume $r,s,k\geq 2$, and $(r,s,k)$ is not $(2,2,2)$. We write $n_{0}=gcd(r+1+s,ks)$. 
Following Saito, we call such $X$ a $k$-sheeted covering of projective space (but note that most varieties covering $\PP^n$ by a degree $k$ finite morphism are not of this form).
\begin{thm}\label{SaitomainA}\cite[Theorem A]{Saio1986}
    If $n_{0}<s$, the $\mathrm{IVHS}$ reconstructs $X$ generically.
\end{thm}

Torelli theorems for some generic (weighted) hypersurfaces are proved using IVHS, and we refer the reader to a series of papers \cite{Don}, \cite{CM_1990__73_2_121_0}, and \cite{Donagi1986}.

\section{Construction of IVHS via Hochschild--Serre algebra of Kuznetsov components}
In this section, first, we reconstruct Infinitesimal Variation Hodge structures of a smooth hypersurface $Y$ of degree $d$ in $\mathbb{P}^n$ such that $(d,n+1)\neq (3,6)$ from its Kuznetsov component.
Then we reconstruct Infinitesimal Variation Hodge structures of generalized Veronese double cones from their Kuznetsov component and then prove the categorical Torelli theorem when they are generic in their moduli space. Finally, we apply a similar method to $k$-sheeted covering of projective spaces and prove the categorical Torelli theorem for them.



\subsection{Smooth hypersurfaces in projective spaces}


Let $Y\subset\mathbb{P}^{n}$ be a degree $d$ hypersurface, and $d< n+1$. Then the semi-orthogonal decomposition of $Y$ is given by 
$$D^b(Y)=\langle\Ku(Y),\oh_Y,\ldots,\oh_Y(n-d)\rangle,$$
where $\Ku(Y)$ is the Kuznetsov component.


\begin{prop}\label{prop_reconstruction_IVHS_quartic_hypersurface}\label{gradealgebrahypersurface}
Let $3\leq d\leq n$, and let $Y,Y'\subset\mathbb{P}^{n}$ be a smooth hypersurface of degree $d$ such that $(d,n+1)\neq (3,6)$, defined by $\omega$ and $\omega'$ respectively. If there is a Fourier--Mukai type equivalence $\Phi:\Ku(Y)\simeq\Ku(Y')$, then 
$\mathrm{IVHS}(Y)\cong\mathrm{IVHS}(Y')$. 
\end{prop}


\begin{proof}
By Proposition~\ref{thm:ordinaryHypersurfacesJacobian}, we have
$$\bigoplus_{t\geq 0}\Jac(\omega)_{td}\cong \bigoplus_{t\geq 0}\Jac(\omega')_{td}.$$
 Thus $\mathrm{IVHS}(Y)\cong \mathrm{IVHS}(Y')$ by Lemma~\ref{weightedIVHS} and Remark~\ref{weightedcommondivisor}.

\end{proof}

\subsection{Generalized Veronese double cones}\label{section_GVDC}
Let $Y\subset\mathbb{P}(\underbrace{1,\ldots,1}_{r+1},2n,kn)$ be the weighted hypersurface of degree $2kn$ such that $k\geq 3$ is odd and $r\geq 2$. Then $Y$ is known as the generalized Veronese double cone in Section~\ref{IVHS_Generalized_double_cone_section}. Note that if $r=2$, $n=1$, and $k=3$, $Y$ is exactly the Veronese double cone. We write $m=\mathrm{gcd}((k+2)n+r+1,2kn)$.

\begin{thm}\label{veroneseIVHS}
  If $m|2n$, the Kuznetsov component $\Ku(Y)$ recovers algebraic data of the infinitesimal variation of the Hodge structure. More precisely, if there is a Fourier--Mukai type equivalence $\Phi:\Ku(Y)\simeq\Ku(Y')$, then $\mathrm{IVHS}(Y)\cong \mathrm{IVHS}(Y')$. 
\end{thm}
\begin{proof}
   We write equation defining $Y$ as $\omega=y^{2}+f$, where $f=x^{k}+\Phi_{k-2}x^{k-2}_{r+2}+\cdots+\Phi_{0}$. Clearly 
   $\mathrm{Jac}(\omega)=\mathrm{Jac}(f)$, and the relations in $\mathrm{Jac}(f)$ has the smallest degree $2kn-2n\geq 2n+1$. According to Theorem \ref{prop_Orlov_identification},
   $$\Ku_{dg}(Y)\simeq \mathrm{Inj}_{\mathrm{coh}}(\mathbb{A}^{r+3},\mathbb{C}^{\ast},\mathcal{O}(2kn),\omega)_{[\underbrace{1,\dots,1}_{r+1},2n,kn]}.$$
   First we show that the Hochschild--Serre algebra of $\Ku(Y)$ reconstructs $\mathrm{Jac}(\omega)_{m}=k[x_{1},x_{2},\cdots,x_{r+1}]_{m}$ and $\mathrm{Jac}(\omega)_{2n}=k[k_{1},x_{2},\cdots,x_{r+1},x]_{2n}$.
Let $\Delta(N)$ be the Fourier--Mukai kernel of the functor
$$-\otimes \mathcal{O}(N):\mathrm{Inj}_{\mathrm{coh}}(\mathbb{A}^{r+3},\mathbb{C}^{\ast},\mathcal{O}(2kn),\omega)_{[\underbrace{1,\dots,1}_{r+1},2n,kn]} \rightarrow \mathrm{Inj}_{\mathrm{coh}}(\mathbb{A}^{r+3},\mathbb{C}^{\ast},\mathcal{O}(2kn),\omega)_{[\underbrace{1,\dots,1}_{r+1},2n,kn]}.$$
The Serre functor $S$ is $-\otimes \mathcal{O}(-(r+1+2n+kn))[r+3]$.
Since $m=gcd(r+1+2n+kn,2kn)$, there exist integers $i'$ and $j'$ such that $$i'(-(r+1+2n+kn))+2knj'=m.$$
We fix such $i'$,$j'$, and take $i_{0}=i',j_{0}=2j'-(r+3)i'$; we then have $S^{i_{0}}[j_{0}]=-\otimes \mathcal{O}(m)$.

By Theorem~\ref{extendedhochschild}, 
\begin{align*}
 \Hom(\Delta,\Delta(N))
 \cong& \bigoplus_{\{g\in \mu_{2kn},\  \operatorname{rk}W_{g}\equiv\ 0\ (mod\ 2)\} }\mathrm{Jac}(\omega_{g})_{N-k_{g}+2kn(\frac{-\operatorname{rk}W_{g}}{2})}
 \end{align*}
We write $\mu_{2kn}=\langle \gamma\rangle$ with $\gamma=\exp{\frac{2\pi i}{2kn}}$. An element $g \in \mu_{2kn}$ is of the form $\gamma^j$ with $0 \le j < 2kn$, and we get the following cases.
\begin{itemize}
\item $j = 0$: $(\mathbb{A}^{r+3})^{\gamma^{j}}=\AA^{r+3}$; $\operatorname{rk} W_{\gamma^{j}}=0$; $k_{\gamma^{j}}=0$, $\Jac(\omega_{\gamma^j}) = \Jac(\omega)$
    \item $2|j$, $k\not| j$: $(\mathbb{A}^{r+3})^{\gamma^{j}}=(0,\cdots,0,y)$; $\operatorname{rk} W_{\gamma^{j}}=r+2$; $k_{\gamma^{j}}=-(r+1+2n)$, $\Jac(\omega_{\gamma^j}) \cong \CC$.

    \item $2\not|j$, $k|j$: $(\mathbb{A}^{r+3})^{\gamma^{j}}=(0,\cdots,0,x,0)$; $\operatorname{rk} W_{\gamma^{j}}=r+2$; $k_{\gamma^{j}}=-(r+1+kn)$, $\Jac(\omega_{\gamma^j}) \cong \CC[x]/(x^{k-1})$.
    \item $2|j$, $k|j$, $j \not= 0$: $(\mathbb{A}^{r+3})^{\gamma^{j}}=(0,\cdots,0,x,y)$; $\operatorname{rk} W_{\gamma^{j}}=r+1$; $k_{\gamma^{j}}=-(r+1)$, $\Jac(\omega_{\gamma^j}) \cong \CC[x]/(x^{k-1})$.
   \item $2\not|j$, $k\not|j$:
   $(\mathbb{A}^{r+3})^{\gamma^{j}}=(0,\cdots,0,0)$; $\operatorname{rk} W_{\gamma^{j}}=r+3$; $k_{\gamma^{j}}=-(r+1)-kn-2n$, $\Jac(\omega_{\gamma^j}) \cong \CC$.
\end{itemize}
Thus, we have for any $N$ that
\[
\Hom(\Delta, \Delta(N)) = \Jac(\omega)_N \oplus \mathbb{C}_{N-((r+1)(kn-1)+(kn-2n))} \oplus \mathbb{C}[x]/(x^{k-1})_{N-(r+1)(kn-1)}.
\]

If $N=m\ \text{or}\ 2n$, then
$$N< (r+1)(kn-1)+(kn-2n)\ \text{and}\ N < (r+1)(kn-1).$$
because $k\geq 3$, $r\geq 2$.
It follows that
\[
\Hom(\Delta,\Delta(m))\cong\mathrm{Jac}(\omega)_{m}
\]
and
\[
\Hom(\Delta, \Delta(2n))\cong\mathrm{Jac}(\omega)_{2n}.
\]

 We show the Hochschild--Serre algebra of $\Ku(Y)$ reconstructs $\bigoplus_{j=0}\mathrm{Jac}(Y)_{jm}$. Define $\mathrm{HS}(Y)_{m}$ as the sub-graded algebra generated by 
   $$\Hom(\mathrm{Id},S^{i_{0}}[j_{0}]),\ \Hom(\mathrm{Id},(S^{i_{0}}[j_{0}])^{\frac{2n}{m}})\subset \bigoplus^{\frac{\sigma}{m}}_{t=0}\Hom(\mathrm{Id},(S^{i_{0}}[j_{0}])^{t}),$$
   where $\sigma$ is the degree of the maximal degree term of $\mathrm{Jac}(\omega)$. 
   Similarly define $\mathrm{HS}(Y')_{m}$.
   Since $\mathrm{Jac}(\omega)_{m}=k[x_{1},x_{2},\cdots,x_{r+1}]_{m}$ and $\mathrm{Jac}(\omega)_{2n}=k[x_{1},x_{2},\cdots,x_{r+1},x]_{2n}$, $\mathrm{Jac}(\omega)_{tm}$ can be generated by $\mathrm{Jac}(\omega)_{m}$ and $\mathrm{Jac}(\omega)_{2n}$. 
   Then we have an isomorphism of graded algebras,
   $$\mathrm{HS}(Y)_{m}\cong \bigoplus^{\frac{\sigma}{m}}_{t=0}\mathrm{Jac}(\omega)_{tm}.$$
   The same results hold for $\mathrm{HS}(Y')_{m}$. Since an equivalence $\Ku(Y)\cong\Ku(Y')$ yields a bigraded isomorphism $\mathrm{HS}(\Ku(Y))\cong \mathrm{HS}(\Ku(Y'))$ by Theorem~\ref{Moritaserrealgebra},
   
   we have an isomorphism of graded algebras,
   $$\bigoplus^{\frac{\sigma}{m}}_{t=0}\mathrm{Jac}(\omega)_{tm}\cong \bigoplus^{\frac{\sigma}{m}}_{t=0}\mathrm{Jac}(\omega')_{tm}.$$
   
   By Lemma \ref{weightedIVHS}, $\bigoplus^{\frac{\sigma}{m}}_{t=0}\mathrm{Jac}(Y)_{tm}$ is enough to recover IVHS of $Y$. Then the desired result follows. 
\end{proof}

\begin{thm}(Generic categorical Torelli)\label{generictorelli}
 Let $Y$ and $Y'$ be generic degree $d$ smooth hypersurface of $\mathbb{P}(1,1,\cdots,1,2n,kn)$ satisfying the assumption in Theorem~\ref{veroneseIVHS}, and $m<2n$. If $\Ku(Y)\cong \Ku(Y')$, then $Y\cong Y'$.
\end{thm}

\begin{proof}
   According to Theorem \ref{veroneseIVHS}, we have an isomorphism of IVHS of $Y$ and $Y'$. By Proposition~ \ref{IVHSgeneric}, IVHS recovers isomorphism classes generically, namely $Y\cong Y'$.  
\end{proof}


\subsection{$k$-sheeted covering over $\mathbb{P}^{r}$}\label{section_kSC}
Let $X$ be a smooth hypersurface of degree $ks$ in $\mathbb{P}(\underbrace{1,1,\cdots,1}_{r+1},s)$ defined by $\omega$. Assume $r,s,k\geq 2$, and $(r,s,k)$ is not $(2,2,2)$. We write $n_{0}=gcd(r+1+s,ks)$.

\begin{thm}\label{thm_generic_Torelli_2}
  If $2\cdot n_{0}<\operatorname
{min}\left( 2s, ((r+2)k-2)s-2(r+1),(ks-2)(r+1)\right)$ and $n_{0}|s$, $\Ku(X)$ reconstructs $X$ generically.  
\end{thm}
\begin{proof}
  Using the same method in Theorem \ref{generictorelli}, if  $$2\cdot gcd(r+1+s, ks)<\operatorname
{min}\left( ((r+2)k-2)s-2(r+1),(ks-2)(r+1)\right),$$
we have $\Hom(\Delta,\Delta(n_{0}))\cong \mathrm{Jac}(\omega)_{n_{0}}$ and $\Hom(\Delta,\Delta(s))\cong \mathrm{Jac}(\omega)_{s}$. Thus the Hochschild--Serre algebra of $\Ku(X)$ reconstructs the graded algebra $\bigoplus^{\frac{\sigma}{n}}_{t=0}\mathrm{Jac}(\omega)_{tn_{0}}$. According to Lemma \ref{weightedIVHS}, the Hochschild--Serre algebra of $\Ku(X)$ reconstructs IVHS of $X$. According to Theorem \ref{SaitomainA}, if we assume $n_{0}<s$ further, the IVHS reconstructs $X$ generically. Thus $\Ku(X)$ reconstructs $X$ generically.
\end{proof}

\section{Categorical Torelli theorems for weighted hypersurfaces}\label{section_improvement}
In the previous section, we reconstruct the infinitesimal Variation Hodge structure (IVHS) for two types of weighted hypersurfaces and prove categorical Torelli theorems for generic ones by Proposition~\ref{IVHSgeneric} and Theorem~\ref{SaitomainA} respectively. In this section, we remove the generic assumption on these weighted hypersurfaces and prove the following theorem:

\begin{thm}\label{thm_categorical_Torelli_thm_weighted_hypersurfaces}
Categorical Torelli theorem holds for two types of weighted hypersurfaces satisfying assumptions in Theorem~\ref{generictorelli} and Theorem~\ref{thm_generic_Torelli_2} respectively. 
\end{thm}

\begin{proof}
We show the categorical Torelli theorem for generalized Veronese double cone in Section~\ref{section_GVDC}. We indeed show the sub-graded algebra $\bigoplus^{\frac{\sigma}{n}}_{t=0}\mathrm{Jac}(\omega)_{tn_{0}}$ is enough to reconstruct the hypersurface. Using notations in Theorem \ref{veroneseIVHS},
$\Jac(\omega)_{m}=\mathbb{C}[x_{0},x_{1},\cdots,x_{r}]_{m}$, and $\Jac(\omega)_{2n}=\langle\mathbb{C}[x_{0},x_{1},\cdots,x_{r}]_{2n},x\rangle$. We have a commutative diagram,
$$\xymatrix{\Jac(\omega)_{m}\times \Jac(\omega)_{m}\ar[r]\ar[d]&\Jac(\omega)_{2m}\ar[d]
\\
\Jac(\omega')_{m}\times \Jac(\omega')_{m}\ar[r]&\Jac(\omega')_{2m}},$$
Where the vertical maps are isomorphisms, and the horizontal maps are multiplication. Since $m<2n$, and $m|2n$, we have $2m\leq 2n$, then $\Jac(\omega)_{2m}$ is $\mathbb{C}[x_{1},x_{2},\cdots,x_{r+1}]_{2m}$ or $\langle\mathbb{C}[x_{1},x_{2},\cdots,x_{r+1}]_{2n},x\rangle$($2m=2n$ in this case). Therefore the polynomial structure of $\mathbb{C}[x_{1},x_{2},\cdots,x_{r+1}]$ can be reconstructed by the commutative diagram above \cite[Lemma 4.2]{Don}. The isomorphism 
$\Jac(\omega)_{2n}\cong \Jac(\omega')_{2n}$ maps $x$ to $c\cdot x+\sum^r_{i=0}a_{i}x^{c_{i_{0}}}_{i}x^{c_{i_{1}}}_{i_{1}}\cdots x^{c_{i_{t}}}_{i_{t}}$($c\neq 0,c_{i_{0}}+c_{i_{1}}++\cdots+c_{i_{t}}=2n$). 
Thus we have an automorphism of $\mathbb{P}(1,1,\cdots,1,2n,kn)$ which induces isomorphism $\Jac(\omega)_{2kn}\cong \Jac(\omega')_{2kn}$. By \cite[Theorem 1.3]{Donagi1986}, $Y\cong Y'$.
 
 A similar method shows the categorical Torelli theorem for $k$-sheeted covering of $\mathbb{P}^r$ in Section~\ref{section_kSC}. 
\end{proof}




\bibliographystyle{alpha}
{\small{\bibliography{Categoricaltorellihypersurface}}}
\end{document}